\documentclass{amsart}

\usepackage{amsmath}
\usepackage{amssymb}

\def\BC{{\mathbb C}}
\def\BD{{\mathbb D}}

\newcommand{\bpr}{{\noindent\textbf{Proof.}\ }}
\newcommand{\epr}{{\hfill $\Box$}}

\newcommand{\sH}{{\mathcal H}}
\newcommand{\sK}{{\mathcal K}}

\newcommand{\sM}{{\mathcal M}}
\newcommand{\sF}{{\mathcal F}}
\newcommand{\sU}{{\mathcal U}}
\newcommand{\sV}{{\mathcal V}}
\newcommand{\sG}{{\mathcal G}}
\newcommand{\sL}{{\mathcal L}}
\newcommand{\sE}{{\mathcal E}}
\newcommand{\sD}{{\mathcal D}}
\newcommand{\sW}{{\mathcal W}}
\newcommand{\sY}{{\mathcal Y}}
\newcommand{\sX}{{\mathcal X}}

\newcommand{\im}{\textup{Im\,}}
\newcommand{\mat}[2]{\ensuremath{\left[\begin{array}{#1}
#2
\end{array} \right]}}

\newcommand{\la}{\lambda}
\newcommand{\om}{\Omega}
\newcommand{\oo}{\omega}
\newcommand{\ga}{\Gamma}

\newcommand{\coup}[1]{\{U_{#1},\tau_{#1}\}}
\newcommand{\coupon}[1]{\{U_{#1}\mbox{ on }\sK_{#1},\tau_{#1}\}}
\newcommand{\lifset}{\{T',A\}}
\newcommand{\lifsetT}{\{T_\circ', A_\circ\}}
\newcommand{\liftset}{\{A, T', U', R, Q\}}
\newcommand{\liftsetV}{\{A, T', V, R, Q\}}
\newcommand{\liftsetVT}{\{A_\circ, T_\circ', \tilde V, R_\circ, Q_\circ\}}

\newcommand{\eL}{{\mathbf{L}}}
\newcommand{\eS}{{\mathbf{S}}}
\newcommand{\bvz}{\bigvee_{n=0}^\infty }

\newtheorem{theorem}{Theorem}[section]
\newtheorem{corollary}[theorem]{Corollary}
\newtheorem{lemma}[theorem]{Lemma}
\newtheorem{proposition}[theorem]{Proposition}

\begin{document}

%
\title[Coupling and relaxed commutant lifting]
 {Coupling and relaxed commutant lifting}

\author[A.E. Frazho]{A.E. Frazho}

\address{%
Department of Aeronautics and Astronautics\\
Purdue University \\
West Lafayette, IN 47907, USA}

\email{frazho@ecn.purdue.edu}

\author{S. ter Horst}
\address{Afdeling Wiskunde,\\
Faculteit der Exacte Wetenschappen,
Vrije Universiteit\\
De Boelelaan 1081a, 1081 HV Amsterdam, The Netherlands}

\email{terhorst@few.vu.nl}

\author{M.A. Kaashoek}
\address{Afdeling Wiskunde\\
Faculteit der Exacte Wetenschappen,
Vrije Universiteit\\
De Boelelaan 1081a, 1081 HV Amsterdam, The Netherlands}

\email{ma.kaashoek@few.vu.nl}



\date{}

\begin{abstract}
A Redheffer type description of the set of all contractive
solutions to the relaxed commutant lifting problem is given. The
description involves a set of Schur class functions which is
obtained by combining the method of isometric coupling with
results on isometric realizations. For a number of special cases,
including the case of the classical commutant lifting theorem, the
description yields a proper parameterization of the set of all
contractive solutions, but examples show that, in general, the
Schur class function determining the contractive lifting does not
have to be unique. Also some  sufficient conditions are given
guaranteeing that the corresponding relaxed commutant lifting
problem has only one solution.
\end{abstract}

\maketitle

\setcounter{section}{-1}
\section{Introduction}\label{sec:intro}
\setcounter{equation}{0}

This paper is devoted to  the relaxed commutant lifting theorem in
\cite{ffk02}. This theorem is a generalization of the classical
commutant lifting theorem \cite{Sz.-NF3}, and it includes as
special cases the Treil-Volberg lifting theorem \cite{tv94}, and
its weighted version due to Biswas, Foias and Frazho \cite{bff99}.

To state the relaxed commutant lifting theorem, let us first
recall the general setup. The starting point is a \emph{lifting
data set} $\liftset$ consisting of five Hilbert space operators.
The operator $A$ is a contraction mapping $\sH$ into ${\sH}'$, the
operator $U'$ on $\sK'$ is a minimal  isometric lifting of $T'$ on
$\sH'$, and $R$ and $Q$ are operators from $\sH_0$ to $\sH$,
satisfying the following constraints
\[
T'AR=AQ  \quad \textup{and}\quad R^*R\leq Q^*Q.
\]
Given this  data set the relaxed commutant lifting theorem in
\cite{ffk02} states that there exists a contraction $B$ from $\sH$
to $\sK'$ such that
\begin{equation}\label{rclt}
\Pi_{\sH'}B=A\quad\mbox{and}\quad U'BR=BQ.
\end{equation}
Here $\Pi_{\sH'}$ is the orthogonal projection from $\sK'$ onto
$\sH'$. In fact,  \cite{ffk02}  provides an explicit construction
for  a contraction  $B$ satisfying  (\ref{rclt}). In the sequel we
say that $B$ is a \emph{contractive interpolant} for $\liftset$ if
$B$ is a contraction from $\sH$ into $\sK$ satisfying
(\ref{rclt}).

In this paper we present a Redheffer type formula to describe the
set of all contractive interpolants for  $\liftset$. In order to
state  our main results we need some auxiliary operators.  To this
end, let $D_\circ$ be the positive square root of $Q^*Q-R^*R$, and
set
\begin{equation} \label{spacesF}
\sF=\overline{D_AQ\sH_0}  \quad \mbox{and}\quad \sF'=
\overline{\left[
\begin{array}{c}
D_\circ\\ D_{T'}AR\\D_AR
\end{array}
\right]\sH_0}.
\end{equation}
Notice that $\sF$ is a subspace of $\sD_A$ and
$\sF'$ is a subspace of $\sD_\circ \oplus \sD_{T'} \oplus \sD_A$.
Here we follow the convention that for a contraction $C$, the
symbol $D_C$ denotes the positive square root of $I - C^*C$ and
$\sD_{C}$ stands for the closure of the range of $D_{C}$.
Furthermore, $\sD_\circ=\overline{D_\circ\sH_0}$. Since $T'AR=AQ$,
we know from formula (4.11) in \cite{ffk02} that there exists a
unique unitary operator $\omega$ mapping $\sF$ onto $\sF'$ such
that
\begin{equation}
\label{omega} \omega(D_AQh)=\left[
\begin{array}{c}
D_\circ\\ D_{T'}AR\\D_AR
\end{array}
\right]h, \quad h\in \sH_0.
\end{equation}
We also need the projections $\Pi_{T'}$ and $\Pi_A$ defined  by
\begin{equation}
\label{projects} \Pi_{T'}=\left[\begin{array}{ccc} 0&I&0
\end{array}\right]: \left[
\begin{array}{c}
\sD_\circ\\ \sD_{T'}\\ \sD_A
\end{array}
\right]\to \sD_{T'}, \quad \Pi_A=\left[\begin{array}{ccc} 0&0&I
\end{array}\right]: \left[
\begin{array}{c}
\sD_\circ\\
\sD_{T'}\\ \sD_A
\end{array}
\right]\to \sD_A.
\end{equation}

Notice that  the previous definitions  only relied upon  the
operators $A,T',R$ and $Q$. The minimal isometric lifting $U'$ did
not play a role. Recall that all minimal isometric liftings of the
same contraction are isomorphic. So without loss of generality, in
our main theorem, we can  assume that $U'=V$ is the
Sz.-Nagy-Sch\"affer minimal isometric lifting of $T'$ which acts
on $\sH'\oplus H^2(\sD_{T'})$. The definitions of a minimal
isometric lifting and the Sz.-Nagy-Sch\"affer lifting are
presented in the next section. Finally, given  Hilbert spaces
$\sU$ and $\sY$, we write $\eS(\sU, \sY)$ for the set of all
operator-valued functions which are analytic on the open unit disk
$\BD$ and whose values are contractions from $\sU$ to $\sY$. We
refer to $\eS(\sU, \sY)$ as the \emph{Schur class} associated with
$\sU$ and $\sY$. We are now ready to state our first main result.

\begin{theorem}
\label{mainth} Let $\liftsetV$  be a lifting data set,  where  $V$
on $\sH'\oplus H^2(\sD_{T'})$ is  the Sz.-Nagy-Sch\"affer minimal
isometric lifting of $T'$. Then  all contractive interpolants for
this date set are given by
\begin{equation}
\label{sols}Bh=\left[
\begin{array}{c}
Ah\\ \Pi_{T'}F(\lambda)(I_{\sD_A}-\lambda \Pi_A F(\lambda))^{-1}D_Ah
\end{array}
\right], \quad h\in \sH,
\end{equation}
where $F$ is any function from the Schur class $\eS(\sD_A,
\sD_\circ\oplus\sD_{T'}\oplus \sD_A)$ satisfying
$F(0)|\sF=\omega$.
\end{theorem}

In general, formula (\ref{sols}) does not establish a one to one
correspondence between $B$ and the parameter $F$. It can happen
that different $F$'s  yield the same $B$. For instance, assume
$\sH_0,\sH$ and $\sH'$ to be equal to $\BC$, let $A,R$ and $Q$ be
the zero operator on $\BC$, and take for $T'$ the identity
operator on $\mathbb{C}$. Since $T'$ is an isometry,  the
Sz.-Nagy-Sch\"affer minimal isometric lifting $V$ of $T'$ is equal
to $T'$. The latter implies that there is only one contractive
interpolant $B$ for the data set $\liftsetV$, namely $B=A$. The
fact that $R$ and $Q$ are the zero operators on $\BC$ implies that
$\sF=\{0\}$ and $\sF'=\{0\}$. It follows that for this data set
$\liftsetV$ the only  contractive interpolant $B$ is given by
formula (\ref{sols}) where for $F$ we can take any function in the
Schur class $\eS(\mathbb{C},\mathbb{C})$. The previous example can
be seen as  a special case of our second main theorem.

\begin{theorem}\label{thmmain2}
Let $B$ be a contractive interpolant for the data set $\liftsetV$
where $V$ is  the Sz.-Nagy-Sch\"affer minimal isometric lifting of
$T'$. Then there is a one to one mapping from the set of all $F$
in $\eS(\sD_A,\sD_\circ\oplus\sD_{T'}\oplus\sD_A)$ with
$F(0)|\sF=\oo$ such that $B$ is given by $(\ref{sols})$ onto the
set $\eS(\sG_B,\sG_B')$, with $\sG_B$ and $\sG_B'$ being given by
\begin{equation}\label{eq:FB0}
\sG_B=\sD_B\ominus\overline{D_BQ\sH_0}\quad\mbox{and}\quad
\sG_B'=(\sD_\circ\oplus\sD_B)\ominus\overline{\mat{c}{D_\circ\\D_BR}\sH_0}.
\end{equation}
\end{theorem}

Our proof of the above theorem also provides a procedure to obtain
a mapping of the type referred to in the theorem.

It is interesting to specify Theorems \ref{mainth} and
\ref{thmmain2} for the case when in the lifting data set
$\liftsetV$ the operators $A$, $T'$, $R$ and $Q$ are zero
operators. In this case the intertwining condition $VBR=BQ$, where
$V$ is the Sz.-Nagy-Sch\"affer minimal isometric lifting of
$T'=0$, is trivially fulfilled, and hence $B$ is a contractive
interpolant if and only if
\[
Bh=\mat{c}{0\\ \Theta(\cdot)h}, \quad h\in\sH,
\]
where $\Theta$ is any function in $H_{ball}^2(\sL(\sH,\sH'))$. The
latter means that $\Theta$ is a $\sL(\sH,\sH')$-valued analytic
 function  on $\mathbb{D}$ such  that for each $h \in \sH$
the function $\Theta(\cdot)h$ belongs to the Hardy space
$H^2(\sH')$, and $\|\Theta(\cdot)h\|_{H^2(\sY)}\leq \|h\|$. It
follows that  Theorems \ref{mainth} and \ref{thmmain2} have the
following corollaries.

\begin{corollary}
\label{firstcor} Let $F$ be any function in the Schur class
$\eS(\sH, \sH'\oplus \sH)$, and let $\Pi$ and $\Pi'$ be the
orthogonal projections of $\sH'\oplus \sH$ on $\sH$ and $\sH'$,
respectively. Then the function $\Theta$ defined by
\begin{equation}
\label{allsols} \Theta(\lambda)={\Pi}'F(\lambda)(I_{\sH}-\lambda
\Pi F(\lambda))^{-1}
\end{equation}
belongs to $H_{ball}^2(\sL(\sH,\sH'))$, and any function in
$H_{ball}^2(\sL(\sH,\sH'))$ is obtained in this way
\end{corollary}

\begin{corollary}
\label{secondcor} Let $\Theta\in H_{ball}^2(\sL(\sH,\sH'))$. Then
there is a one to one mapping from the set of all $F$ in $\eS(\sH,
\sH'\oplus \sH)$ such that $(\ref{allsols})$ holds onto the set
$\eS(\sD_\Gamma, \sD_\Gamma)$, where $\Gamma$ is the contraction
from $\sH$ into $H^2(\sH')$ defined by
\[
(\Gamma h)(\lambda)  = \Theta(\lambda) h, \qquad h\in \sH, \
\lambda \in \mathbb{D}.
\]
\end{corollary}

When $\sH=\sH'=\mathbb{C}$, and hence $\Theta$ is a scalar
function, Corollary \ref{firstcor} can be found in \cite{Sar89},
page 490, provided $\Theta$ is of unit $H^2$ norm. For $p\times q$
matrix functions $\Theta$, when $\sH= \mathbb{C}^q$ and $\sH'=
\mathbb{C}^p$,  Corollary \ref{firstcor} is Theorem 2.2 in
\cite{ABP95}. For the general operator valued case Corollary
\ref{firstcor} seems to be new. Corollary \ref{secondcor} seems to
be new even in the scalar case. Notice that in the scalar case the
space $\sD_\Gamma$ in Corollary \ref{secondcor} consists of the
zero element only if $\Theta$ is of unit $H^2$ norm, and
$\sD_\Gamma=\mathbb{C}$ otherwise.

Another case of special interest is the classical commutant
lifting problem. As we know from \cite{ffk02} the commutant
lifting theorem can be obtained by applying the relaxed commutant
lifting theorem to the data set $\{A,T',U',I_\sH,Q\}$ where
$\sH_0=\sH$, the operator $R$ is the identity operator on $\sH$
and $Q$ is an isometry; see \cite{ffk02}. In this case,  the space
$\sG_B'$ in Theorem \ref{thmmain2} consists of the zero element
for any choice of the contractive interpolant $B$. In other words,
for the case of the classical commutant lifting formula
(\ref{sols}) provides a \emph{proper parameterization}, that is,
for every contractive interpolant $B$ for $\liftsetV$ there exists
a unique $F$ in $\eS(\sD_A,\sD_\circ\oplus\sD_{T'}\oplus\sD_A)$
with $F(0)|\sF=\oo$ such that $B$ is given by (\ref{sols}).
Finally, it is noted that this formula also yields  the Redheffer
type parameterization
 for the commutant lifting theorem presented
in Section  XIV of \cite{ff}.

If in Theorem \ref{mainth} we take $F(\lambda)\equiv
\omega\Pi_\sF$, where $\Pi_\sF$ is the orthogonal projection of
$\sD_A$ onto $\sF$, then the contractive interpolant $B$ in
(\ref{sols})  is precisely the central solution presented in
\cite{ffk02}.

{}From Theorem \ref{mainth} we see that $\sF=\sD_A$ implies that
there is a unique contractive interpolant  (which is known from
Theorem 3.1 in \cite{ffk02}). Other conditions of uniqueness will
be given in the final section of the paper.

We shall prove Theorems \ref{mainth} and \ref{thmmain2} by
combining the method of isometric coupling with some aspects of
isometric realization theory. The theory of isometric couplings
originates from \cite{AA65}, \cite{AA66}, and was used to study
the commutant lifting problem for the first time in \cite{ar83a}
-- \cite{ar83e}; see also, Section VII.7 in \cite{ff}.

The paper consist of six sections not counting this introduction.
The first two sections have a preliminary character, and review
the notions of an isometric lifting (Section \ref{sec:IL}), and an
isometric realization (Section \ref{sec:IR}). In the third section
we develop the notion of an isometric coupling of a pair of
contractions which provides  the main tool in this paper. In
Section \ref{sec:RCLT0} we prove Theorem \ref{mainth} for the case
when $R^*R=Q^*Q$, and in Section \ref{sec:RCLT1} we prove Theorem
\ref{mainth} in its full generality. In the final section we prove
Theorem \ref{thmmain2}, and we present a few sufficient conditions
for the case when (\ref{sols}) provides a proper parameterization,
and also conditions for uniqueness of the solution.

We conclude this introduction with a few words about notation and
terminology. Throughout capital calligraphic letters denote
Hilbert spaces. The Hilbert space direct sum of $\sU$ and $\sY$ is
denoted by
\[
\sU\oplus \sY \quad \mbox{or by}\quad \left[
\begin{array}{c}
\sU\\ \sY
\end{array}
\right].
\]
The set of all bounded
linear operators from $\sH$ to $\sH'$ is denoted by $\eL(\sH,
\sH')$. The identity operator on the space $\sH$ is denoted by
$I_{\sH}$ or just by $I$, when the underlying space is clear from
the context. By definition, a \emph{subspace} is a closed linear
manifold. If $\sM$ is a subspace of $\sH$, then $\sH\ominus \sM$
stands for the orthogonal complement of $\sM$ in $\sH$. Given a
subspace $\sM$ of $\sH$, the symbol $\Pi_{\sM}$ will denote the
orthogonal projection of $\sH$ onto $\sM$ viewed as an operator
from $\sH$ to $\sM$, and $P_{\sM}$ will denote the orthogonal
projection of $\sH$ onto $\sM$ viewed as  an operator on $\sH$.
Note that $\Pi_{\sM}^*$ is the canonical embedding from $\sM$ into
$\sH$, and hence $P_{\sM}=\Pi_{\sM}^*\Pi_{\sM}$. Instead of
$\Pi_{\sM}^*$ we shall sometimes write $E_{\sM}$, where the
capital $E$ refers to embedding. A subspace $\sM$ of $\sH$ is said
to be \emph{cyclic} for an operator $T$ on $\sH$ whenever
\[
\sH=\bigvee_{n=0}^\infty
T^n\sM=\overline{\textup{span}\,\{T^n\sM\mid n=0, 1,2, \ldots\}}.
\]
Finally, by definition, a $\eL(\sH, \sH')$-valued \emph{Schur
class function} is a function in $\eS(\sH,\sH')$, i.e., an
operator-valued function which is analytic on the open unit disk
$\BD$ and whose values are contractions from $\sH$ to $\sH'$.

\section
{Isometric liftings}\label{sec:IL} \setcounter{equation}{0} In
this section we review  some facts concerning isometric liftings
that are used throughout this paper. For a more complete account
we refer to the book \cite{sznf} (see also Chapter VI in
\cite{ff}, and Section 11.3 in \cite{ffk02}).

Let $T'$ on $\sH'$ be a contraction. Recall that an operator $U$
on $\sK$ is a \emph{isometric lifting} of $T'$ if $\sH'$ is a
subspace of $\sK$ and $U$ is an isometry satisfying $\Pi_{\sH'} U
= T^\prime\Pi_{\sH'}$. Isometric liftings exist. In fact, the {\em
Sz.-Nagy-Sch\"{a}ffer   isometric lifting} $V$ of $T^\prime$ is
given by
\begin{equation}\label{eq:NSL}
V = \left[
\begin{array}{cc}
T^\prime             &     0   \\
E D_{T^\prime}         &     S
\end{array} \right] \mbox{ on }
\left[ \begin{array}{c} \sH^\prime \\ H^2(\sD_{T^\prime})
\end{array} \right].
\end{equation}
Here $S$ is the unilateral shift on the Hardy space
$H^2(\sD_{T^\prime})$ and $E$ is the canonical embedding of
$\sD_{T^\prime}$ onto the space of constant functions in
$H^2(\sD_{T^\prime})$. To see that $V$ in (\ref{eq:NSL}) is an
isometric lifting of $T'$  note that any operator $U$ on
$\sK=\sH'\oplus\sM$ is an isometric lifting of $T'$ if and only if
$U$ admits an operator matrix representation of the form
\begin{equation}
\label{eq:U1} U = \left[\begin{array}{cc}
  T' & 0 \\
  Y_1D_{T'} & Y_2
\end{array}\right]\mbox{ on }\left[\begin{array}{c}
  \sH' \\
  \sM
\end{array}\right]\mbox{ where } Y=\mat{cc}{Y_1&Y_2}: \left[\begin{array}{c}
 \sD_{T'} \\
  \sM
\end{array}\right] \to \sM
\end{equation}
is an isometry.

An isometric lifting $U$ of $T'$ is called \emph{minimal} when
$\sH'$ is cyclic for $U$. The Sz.-Nagy-Sch\"{a}ffer isometric
lifting of $T'$ is minimal. If the isometric lifting $U$ is given
by (\ref{eq:U1}), then the lifting is minimal if and only if  the
space $Y_1\sD_{T'}$ is cyclic for $Y_2$.

Two isometric liftings $U_1$ on $\sK_1$ and $U_2$ on $\sK_2$ of
$T'$ are said to be \emph{isomorphic} if  there exists a unitary
operator $\Phi$ from $\sK_1$ onto $\sK_2$ such that
\[
\Phi U_1=U_2\Phi\quad\mbox{and}\quad\Phi h=h\mbox{ for all
}h\in\sH'.
\]
Minimality of an isometric lifting is preserved under an
isomorphism, and two minimal isometric liftings of $T'$ are
isomorphic.

Finally, when $U$ on $\sK$ is a isometric lifting of $T'$, then
the subspace $\sK'$, given by
\[
\sK'=\bigvee_{n=0}^\infty U^n\sH',
\]
is reducing for $U$, that is, both $\sK'$ and its  orthogonal
complement $\tilde{\sK}=\sK\ominus\sK'$ are invariant under $U$.
Furthermore, in that case the operator $U'=\Pi_{\sK'}U|\sK'$ on
$\sK'$ is a minimal isometric lifting of $T'$, and the operator
$U$ admits a operator matrix decomposition of the form
\begin{equation}\label{eq:U2}
U=\mat{cc}{U'&0\\0&\tilde{U}} \mbox{ on }\mat{c}{\sK'\\ \tilde{\sK}},
\end{equation}
where $\tilde{U}$ is an isometry on $\tilde{\sK}$. We shall call
$U'$ in (\ref{eq:U2}) the \emph{minimal isometric lifting of} $T'$
\emph{associated with} $U$.

The following proposition summarizes the results referred to above
in a form that will be convenient for this paper. For details we
refer to  Section 11.3 in \cite{ffk02}.

\begin{theorem}\label{thmlift}
Let $T'$ be a contraction on $\sH'$, let $V$ on $\sH'\oplus
H^2(\sD_{T'})$ be the Sz.-Nagy-Sch\"affer (minimal) isometric
lifting of $T'$, and let $U$ on $\sH\oplus \sM$ be an arbitrary
isometric lifting of $T'$ given by \textup{(\ref{eq:U1})}. Then
there exists a unique isometry $\Phi$ from $\sH'\oplus
H^2(\sD_{T'})$ into $\sH'\oplus\sM$ such that $U\Phi=\Phi V$ and
$\Phi|\sH'=I_{\sH'}$. In fact, $\Phi$ is given by
\[
\Phi=\mat{cc}{I_{\sH'}&0\\0&\Lambda} :\mat{c}{\sH'\\
H^2(\sD_{T'})}\to\mat{c}{\sH'\\ \sM},
\]
where $\Lambda$ is defined by
\[
\Lambda h=\sum_{n=0}^\infty Y^n_2Y_1h_n,\quad
h(\la)=\sum_{n=0}^\infty \la^nh_n\in H^2(\sD_{T'}),
\]
with $Y_1$ and $Y_2$ as in~$(\ref{eq:U1})$. Moreover,
$(\Lambda^*m)(\la)=Y_1^*(I-\la Y_2^*)^{-1}m$ for each $m\in\sM$.
Finally, $\Phi$ is unitary if and only if $U$ is a minimal
isometric lifting of $T'$, and in that case the isometric liftings
$V$ and $U$ of $T'$ are isomorphic.
\end{theorem}

The isometry  $\Phi$ introduced in the above theorem will be
referred to as the \emph{unique isometry associated with} $T'$
\emph{that intertwines} $V$ \emph{with} $U$. Since $V$ is uniquely
determined by $T'$, we shall denote this isometry simply by
$\Phi_{U,\,T'}$.  When $U$ on $\sK$ is an isometric lifting of
$T'$ and $U'$ on $\sK'$ is the minimal isometric lifting  of $T'$
associated with $U$, then the operator $\Pi_{\sK'}\Phi_{U,\,T'}$
is the unique isometry associated with $T'$ that intertwines $V$
with $U'$, that is, $\Phi_{U',\,T'}=\Pi_{\sK'}\Phi_{U,\,T'}$ or,
equivalently, $\Pi^*_{\sK'}\Phi_{U',\,T'}=\Phi_{U,\,T'}$.

\section{Isometric realizations}\label{sec:IR}
\setcounter{equation}{0}

In this section we review some of the classical results on
controllable isometric realizations, and we prove a few additional
results that will be useful in the later sections.

We say that $\{Z,B,C,D; \sX,\sU,\sY\}$ (or simply $\{Z,B,C,D\}$)
is a {\em realization} of a $\eL(\sU,\sY)$-valued function $G$ if
\begin{equation}\label{eq:IR1}
G(\la)=D+\la C(I_\sX-\la Z)^{-1}B
\end{equation}
for all $\la$ in some open neighborhood of the origin in the
complex plane. Here  $Z$ is an operator on $\sX$  and $B$ is an
operator from $\sU$ into $\sX$ while $C$ is an operator mapping
$\sX$ into $\sY$ and $D$ is an operator from  $\sU$ into $\sY$
(where $\sX$, $\sU$ and $\sY$ are all Hilbert spaces). In this
case, we refer to the function defined by the right hand side of
(\ref{eq:IR1}) as the associated \emph{transfer function}. A
realization $\{Z,B,C,D\}$ is called \emph{isometric} if the
operator
\begin{equation}\label{Mabcd0}
M = \left[
\begin{array}{cc}
  D & C \\
  B & Z \\
\end{array}\right]:
\left[\begin{array}{c}
  \sU \\
  \sX
\end{array}\right]\rightarrow
\left[\begin{array}{c}
  \sY \\
  \sX
\end{array}\right]
\end{equation}
is an isometry. The $2\times 2$ operator matrix in (\ref{Mabcd0})
is called the \emph{system matrix associated with the realization}
$\{Z,B,C,D\}$. The transfer function of an isometric realization
belongs to the Schur class $\eS(\sU,\sY)$, that is, if
$\{Z,B,C,D\}$ is an isometric realization, then the function $G$
defined by (\ref{eq:IR1}) is a contractive analytic
$\eL(\sU,\sY)$-valued function on $\BD$. Conversely, if $G\in
\eS(\sU,\sY)$, then there is an isometric realization $\{Z
,B,C,D\}$ such that (\ref{eq:IR1}) holds for all $\lambda$ in
$\BD$.

The transfer function of a realization can also be expressed in
terms of the system matrix $M$. In fact, if $\{Z,B,C,D\}$ is a
realization and $M$ is the associated system matrix, then in a
neighborhood of the origin the transfer function $G$ is also given
 by
\begin{equation}\label{eq:CAF-ISM}
G(\la)=\Pi_\sY M(I_{\sU\oplus\sX}-\lambda J_{\sX}M)^{-1}\Pi_\sU^*,
\end{equation}
where $J_\sX$ is the partial isometry from $\sY\oplus\sX$ to
$\sU\oplus\sX$ given by
\[
J_\sX=\mat{cc}{0&0\\0&I_\sX}:\mat{c}{\sY\\\sX}\to\mat{c}{\sU\\\sX}.
\]
Indeed, for $\lambda$ sufficiently close to zero we have
\begin{eqnarray*}
G(\la)&=& D+\la C(I_\sX-\la Z)^{-1}B= \mat{cc}{D&C}\mat{c}{I_\sU\\\la(I_\sX-\la Z)^{-1}B}\\
&=&\mat{cc}{D&C}\mat{cc}{I_\sU&0\\\la(I_\sX-\la
Z)^{-1}B&(I_\sX-\la Z)^{-1}}
\mat{c}{I_\sU\\0}\\
&=&\mat{cc}{D&C}\mat{cc}{I_\sU&0\\ -\la B&I_\sX-\la Z}^{-1}\mat{c}{I_\sU\\0}\\
&=&\mat{cc}{D&C} \left(
I_{\sU\oplus\sX}-\la\mat{cc}{0&0\\0&I_\sX}\mat{cc}{D&C\\
B&Z}\right)^{-1}
\mat{c}{I_\sU\\0}\\
&=& \Pi_\sY M(I_{\sU\oplus\sX}-\la J_\sX M)^{-1}\Pi_\sU^* \,.
\end{eqnarray*}
Since the right side of (\ref{eq:IR1}) is a Schur class function
if $M$ in (\ref{Mabcd0}) is an isometry, the same holds true for
the right hand side of (\ref{eq:CAF-ISM}). Notice that the
function $G$ defined by (\ref{eq:CAF-ISM}) can also be written in
the form:
\[
G(\la)=\Pi_\sY (I_{\sY\oplus\sX}-\lambda M
J_{\sX})^{-1}M\Pi_\sU^*.
\]

If for a realization $\{Z ,B,C,D; \sX,\sU,\sY\}$ the space $
\bigvee_{n=0}^\infty Z^nB\sU$ is equal to $\sX$, then the
realization or the pair $\{Z,B\}$ is called \emph{controllable}.
In other words, a realization is controllable if and only if the
space $\overline{B\sU}$ is cyclic for $Z$. In terms of the system
matrix $M$ in (\ref{Mabcd0}) the realization $\{Z ,B,C,D\}$ is
controllable if and only if
\begin{equation}\label{eq:contr2}
\sX=\Pi_\sX\bigvee_{n=0}^\infty (J_\sX M)^n \left[\begin{array}{c}
\sU\\ \{0\}
\end{array}\right].
\end{equation}
The above condition (\ref{eq:contr2}) is also equivalent to the
requirement that $\{J_\sX M,\Pi^*_\sU\}$ is a controllable pair.
In the particular case when $\sU=\sY$ in (\ref{Mabcd0}), condition
(\ref{eq:contr2}) can be written in an even simpler form. This is
the contents of the next lemma.
\begin{lemma}
\label{lem:cont}Let $M$ be as in \textup{(\ref{Mabcd0})}, and
assume $\sU=\sY$. Then $\{Z,B\}$ is controllable if and only if
$\sU \oplus \{0\}$ is cyclic for $M$, that is,
\begin{equation}\label{contr3}
\bigvee_{n=0}^\infty M^n \left[\begin{array}{c} \sU\\ \{0\}
\end{array}\right]=\left[\begin{array}{c}
\sU\\ \sX
\end{array}\right].
\end{equation}
\end{lemma}
\bpr Let $E_\sU$ be the canonical embedding of $\sU$ into
$\sU\oplus\sX$, and define $M_0$ to be the operator
\[
M_0=\left[\begin{array} {cc} 0&0\\ B&Z
\end{array}\right]:\left[\begin{array} {c} \sU \\ \sX
\end{array}\right]\to \left[\begin{array} {c} \sU \\ \sX
\end{array}\right].
\]
Then $M_0=M-E_\sU\left[\begin{array} {cc} C&D
\end{array}\right]$. This feedback relation implies that the pair
$\{M_0,E_\sU\}$ is controllable if and only if the pair
$\{M,E_\sU\}$ is controllable. Thus (\ref{contr3}) holds if and
only if $\{M_0,E_\sU\}$ is controllable. Now notice that
 for all integers $n\geq 1$, we have
\[
 M_0^nE_\sU =\left[\begin{array} {cc} 0&0\\ Z^{n-1}B&Z^n
\end{array}\right]\left[\begin{array} {c} I_\sU \\ 0
\end{array}\right]= \left[\begin{array} {c} 0 \\ Z^{n-1}B
\end{array}\right]:\sU\to \left[\begin{array} {c} \sU \\ \sX
\end{array}\right].
\]
It follows that
\[
\bigvee_{n=0}^\infty M_0^nE_\sU \sU=\left[\begin{array} {c} \sU \\
\{0\}\end{array}\right]\oplus\bigvee_{n=1}^\infty
M_0^nE_\sU=\sU\oplus\bigvee_{n=1}^\infty Z^{n-1}B\sU.
\]
We conclude that (\ref{contr3}) holds if and only if the pair
$\{Z,B\}$ is controllable.\epr

\medskip
A realization $\{Z ,B,C,D\}$ or the pair $\{C,Z\}$ is called
\emph{observable} if $CZ^n x=0$ for all integers $n \geq 0$
implies that the vector  $x$ is equal to zero. Since the
orthogonal complement of $\textup{Ker}\,CZ^n$ is equal to the
closure of $\im (Z^*)^nC^*$, we see that   observability of the
realization $\{Z ,B,C,D\}$ is equivalent to the controllability of
the \emph{dual realization} $\{Z^* ,C^*,B^*,D^*\}$.

Two  realizations $\{Z_1\mbox{ on } \sX_1,B_1,C_1,D_1\}$ and
$\{Z_2\mbox{ on } \sX_2,B_2,C_2,D_2\}$  are said to be {\em
unitarily equivalent} if $D_1 = D_2$ and there exists a unitary
operator $W$ mapping $\sX_1$ onto $\sX_2$ such that
\[
W Z_1 = Z_2 W,\quad W B_1 = B_2 \quad \mbox{and}\quad C_2 W = C_1.
\]
Unitary equivalence does not change the transfer function. More
precisely, when two realizations are unitary equivalent, then
their transfer functions coincide in a neighborhood of zero. For
isometric controllable realizations the converse is also true. In
fact we have the following theorem.

\begin{theorem}\label{th:CAF-IR}
Let $G$ be a  $\eL(\sU,\sY)$-valued function. Then
$G\in\eS(\sU,\sY)$ if and only if $G$ admits an isometric
realization. In this case, $G$ admits a controllable isometric
realization and all controllable isometric realizations of $G$ are
unitarily equivalent. In particular, formula $(\ref{eq:IR1})$
provides a one to one correspondence between  the
$\eL(\sU,\sY)$-contractive analytic functions on $\BD$ and (up to
unitary equivalence) the controllable isometric realizations of
$\eL(\sU,\sY)$-valued functions.
\end{theorem}

The above result appears in a somewhat different form in
\cite{Sz.-NF3} as a theorem representing a Schur class function as
a characteristic operator function. A full proof, with isometric
systems replaced by their dual ones, can be found in \cite{Ando}
which also gives additional references. In Section 1.3 of
\cite{fk03a} the theorem is proved using the Naimark dilation
theory.

\medskip
We conclude this section with a proposition that will be useful in
the later sections. The starting point is an isometry $Y$ of the
type appearing in (\ref{eq:U1}). More precisely,
\begin{equation}
\label{isoY} Y=\mat{cc}{Y_1&Y_2}: \left[\begin{array}{c}
 \sD' \\
  \sM
\end{array}\right] \to \sM.
\end{equation}

\begin{proposition}
\label{propFY}Let $Y$ in $(\ref{isoY})$ be an isometry. Assume
$\sM=\sD\oplus\sX$, and let $\Pi_\sD$ and $\Pi_\sX$ be the
orthogonal projections of $\sM$ onto $\sD$ and $\sX$,
respectively. Put
\begin{equation}
\label{eq:YF}F(\la)=\Pi_{\sD'\oplus\sD}Y^*(I_\sM-\la J_\sX^\prime
Y^*)^{-1}\Pi_\sD^*, \quad \la\in \BD,
\end{equation}
where $\Pi_{\sD'\oplus\sD}$ is the orthogonal projection of
$\sD'\oplus\sM$ onto $\sD'\oplus\sD$,  and
\[J_\sX^\prime: \sD'\oplus\sM\to \sM, \quad
J_\sX^\prime(d'\oplus m)=\Pi_\sX m.
\]
Then $F$ belongs to the  Schur class $\eS(\sD, \sD'\oplus\sD)$ and
\begin{equation}
\label{fundeq*}Y_1^*(I_\sM-\la
Y_2^*)^{-1}\Pi_\sD^*=\Pi'F(\la)\Big(I_\sD-\la\Pi F(\la)\Big)^{-1},
\quad \la\in \BD,
\end{equation}
where $\Pi$ and $\Pi'$ are the orthogonal projections of
$\sD'\oplus\sD$ onto $\sD$ and $\sD'$, respectively.
\end{proposition}

It will be convenient first to prove a lemma. Let $\ga$ be a
contraction from $\sM$ into $\sE_1\oplus\sM$. Partition $\ga$ as a
$2\times1$ operator matrix, as follows
\begin{equation}
\label{matga1} \ga=\left[\begin{array} {c} \ga_1\\ \ga_2
\end{array}\right]:\sM\to\left[\begin{array} {c} \sE_1\\ \sM
\end{array}\right].
\end{equation}
Furthermore, let $\sE_2$ be a subspace of $\sM$, and consider the
function
\begin{equation}
\label{eq:LAM}
\Xi(\la)=\ga_1(I_{\sM}-\la\ga_2)^{-1}\Pi_{\sE_2}^*,\quad
\la\in\BD,
\end{equation}
Here $\Pi_{\sE_2}$ is the orthogonal projection of $\sM$ onto
$\sE_2$. Since $\ga$ is a contraction, the same holds true for
$\ga_2$, and hence $I-\la\ga_2$ is invertible for each
$\la\in\BD$. Thus $\Xi$ is well-defined on $\BD$. Next, let $\sX$
be the orthogonal complement of $\sE_2$ in $\sM$, and thus
$\sM=\sE_2\oplus\sX$. Then $\ga$ also admits a $3\times2$ operator
matrix representation, namely
\begin{equation}
\label{matga2} \ga=\left[\begin{array} {cc} D_1&C_1\\ D_2&C_2\\
B&Z
\end{array}\right]:\left[\begin{array} {c} \sE_2\\ \sX
\end{array}\right]\to\left[\begin{array} {c} \sE_1\\ \sE_2\\ \sX
\end{array}\right].
\end{equation}
Put
\begin{equation}
\label{eq:functF} F(\la)=\left[\begin{array} {c} D_1\\ D_2
\end{array}\right]+\la\left[\begin{array} {c} C_1\\ C_2
\end{array}\right](I_\sX-\la Z)^{-1}B,\quad \la\in\BD.
\end{equation}
Again, since $\ga$ is a contraction, the operator $Z$ is a
contraction, and hence $F$ is well-defined on $\BD$.

\begin{lemma}
\label{lemLAMF} Let $\Xi$ and $F$ be the functions defined by
$(\ref{eq:LAM})$ and $(\ref{eq:functF})$, respectively. Then $F$
belongs to the Schur class $\eS (\sE_2,\sE_1\oplus\sE_2)$ and
\begin{equation}
\label{fundident} \Xi(\la)=\Pi_1F(\la)(I-\Pi_2F(\la))^{-1}\quad
\la\in\BD,
\end{equation}
where $\Pi_1$ and $\Pi_2$ are the orthogonal projections of
$\sE_1\oplus\sE_2$ onto $\sE_1$ and $\sE_2$, respectively.
\end{lemma}
\bpr The function $F$ is the transfer function of the system
\[
\big\{Z,B,\left[\begin{array} {c} C_1\\ C_2
\end{array}\right],\left[\begin{array} {c} D_1\\ D_2
\end{array}\right];\sX,\sE_2,\sE_1\oplus\sE_2\big\}.
\]
By (\ref{matga2}) the system matrix corresponding to this system
is equal to $\ga$, and hence it is a contraction. This implies
that $F$ belongs to the Schur class $\eS(\sE_2,
\sE_1\oplus\sE_2)$; cf., Theorem 4.1 in \cite{akp98} where this is
proved for time-variant systems.

To prove (\ref{fundident}) fix $\la\in\BD$. Using the partitioning
of $\ga$ in (\ref{matga2}) we see that for each $e\in\sE_2$ we
have
\[
\Xi(\la)e=\left[\begin{array} {cc} D_1&C_1
\end{array}\right]\bigg(I-\la \left[\begin{array} {cc}
D_2&C_2\\B&Z \end{array}\right]\bigg)^{-1}\left[\begin{array} {c}
e\\ 0
\end{array}\right].
\]
To find $\Xi(\la)e$ we have to compute the first column of the
inverse of the $2\times2$ operator matrix
\begin{equation}
\label{eq:opmatrix} \left[\begin{array} {cc} I-\la D_2& -\la C_2\\
-\la B& I-\la Z \end{array}\right].
\end{equation}
Since $I-\la Z$ is invertible, the Schur complement $\Delta(\la)$
of $I-\la Z$ in (\ref{eq:opmatrix}) is well-defined and is given
by
\[
\Delta(\la):=I - \la D_2-\la^2C_2(I-\la Z)^{-1}B=I-\la\Pi_2F(\la).
\]
It follows (cf., Remark 1.2 in \cite{bgk}) that
\[\left[\begin{array} {cc} I-\la D_2& -\la C_2\\
-\la B& I-\la Z \end{array}\right]^{-1}=
\left[\begin{array} {cc} (I-\la\Pi_2F(\la))^{-1}& *\\
\la(I-\la Z)^{-1}B(I-\la\Pi_2F(\la))^{-1}& *
\end{array}\right].
\]
Thus
\begin{eqnarray*}
\Xi(\la)e&=&\left[\begin{array} {cc} D_1&C_1
\end{array}\right]\left[\begin{array} {c} (I-\la\Pi_2F(\la))^{-1}e\\
\la(I-\la Z)^{-1}B(I-\la\Pi_2F(\la))^{-1}e
\end{array}\right]\\
\noalign{\vskip6pt} &=& (D_1-\la C_1(I-\la
Z)^{-1}B)(I-\la\Pi_2F(\la))^{-1}e\\
\noalign{\vskip6pt} &=& \Pi_1F(\la)(I-\la\Pi_2F(\la))^{-1}e.
\end{eqnarray*}
Since $e$ is an arbitrary element of $\sE_2$, this proves
(\ref{fundident}). \epr

\medskip\noindent \textbf{Proof of Proposition \ref{propFY}.} Since
$Y$ is assumed to be an isometry, $Y^*$ is a contraction. Now
apply Lemma \ref{lemLAMF} with $\sD'$ in place of $\sE_1$, with
$Y^*$ in place of the contraction $\Gamma$ in (\ref{matga1}), and
with $\sD$ in place of $\sE_2$. With these choices the function
$\Xi$ in (\ref{eq:LAM}) coincides with the function defined by the
left hand side of (\ref{fundeq*}). Thus in order to finish the
proof it remains to show that with $\Gamma=Y^*$, $\sE_1=\sD'$, and
$\sE_2 = \sD$ the function $F$ in (\ref{eq:functF}) is also given
by (\ref{eq:YF}). But this follows by applying to $F$ in place of
$G$ that the function $G$ in (\ref{eq:IR1}) is also given by
(\ref{eq:CAF-ISM}). Indeed, since $F$ is the transfer function of
the system
\[
\big\{Z,B,\left[\begin{array} {c} C_1\\ C_2
\end{array}\right],\left[\begin{array} {c} D_1\\ D_2
\end{array}\right];\sX,\sD,\sD'\oplus\sD\big\},
\]
and the system matrix of this system is equal to $Y^*$, the
equivalence between (\ref{eq:IR1}) and (\ref{eq:CAF-ISM}) yields
in a straightforward way that $F$ in (\ref{eq:functF}) is also
given by (\ref{eq:YF}).\epr

\section{Isometric couplings}\label{sec:IC}
\setcounter{equation}{0} Throughout this section $\lifset$ is a
pair of contractions, $T'$ on a Hilbert space $\sH'$ and $A$ from
a Hilbert space $\sH$ to $\sH'$.

An {\em isometric coupling of} $\lifset$ is a pair $\coupon{}$ of
operators such that $U$ is an isometric lifting of $T'$, acting on
$\sK$ (and thus $\sH'\subset \sK$), and $\tau$ is an isometry from
$\sH$ to $\sK$ with $\Pi_{\sH'}\tau=A$. If the space $\sK$ is of
no interest, then we will just write $\coup{}$. An isometric
coupling $\coupon{}$ of $\lifset$ is called \emph{minimal} if, in
addition, the space $\sH'\vee\tau\sH$ is cyclic for $U$, that is,
\[
\sK=\bigvee_{n=0}^\infty U^n(\sH'\vee\tau\sH).
\]

There exist minimal isometric couplings of $\lifset$. To see this,
let $U$ be the operator on $\sH'\oplus H^2(\sD_{T'})\oplus
H^2(\sD_A)$ given by the following operator matrix representation
\[
U=\mat{ccc}{T'&0&0\\ E_{\sD_{T'}}&S_{\sD_{T'}}&0\\0&0&S_{\sD_A}}.
\]
Here $E_{\sD_{T'}}$ is the canonical embedding of $\sD_{T'}$ onto
the space of constant functions of $H^2(\sD_{T'})$, and
$S_{\sD_{T'}}$ and $S_{\sD_A}$ are the unilateral shifts on
$H^2(\sD_{T'})$ and $H^2(\sD_A)$, respectively. Notice that the
operator defined by the $2\times 2$ operator matrix in the left
upper corner of the matrix for $U$ is the Sz.-Nagy-Sch\"affer
minimal isometric lifting of $T'$. Since $S_{\sD_A}$ is an
isometry, we conclude that $U$ is  also  an isometric lifting of
$T'$. Now let $\tau$ be the isometry defined by
\[
\tau=\mat{c}{A\\0\\ E_{\sD_A}D_A} :\sH\to\mat{c}{\sH'\\
H^2(\sD_{T'})\\ H^2(\sD_A)}
\]
where $E_{\sD_A}$ is the canonical embedding of $\sD_A$ onto the
space of constant functions of $H^2(\sD_A)$. Then $\coup{}$ is a
minimal isometric coupling of $\lifset$.

Two isometric couplings $\coupon{1}$ and $\coupon{2}$ of $\lifset$
are said to be \emph{isomorphic} if there exists a unitary
operator $\Psi$ from $\sK_1$ to $\sK_2$ such that
\[
\Psi U_1=U_2\Psi,\quad\Psi\tau_1=\tau_2\quad\mbox{and} \quad\Psi
h=h\mbox{ for all }h\in\sH'.
\]
In this case
\begin{equation}
\label{eq:morph3a}\Psi U_1\tau_1 = U_2\tau_2.
\end{equation}

Minimality is preserved under isomorphic equivalence. Indeed, when
the pairs  $\coupon{1}$ and $\coupon{2}$ are isomorphic isometric
couplings of $\lifset$, and $\Psi$ from $\sK_1$ to $\sK_2$ is an
isomorphism between the two isometric couplings, then
\begin{eqnarray*}
\bigvee_{n=0}^\infty U_2^n(\sH'\vee\tau_2\sH) &=&
\bigvee_{n=0}^\infty (\Psi U_1\Psi^*)^n(\sH'\vee\Psi\tau_1\sH)=
\bigvee_{n=0}^\infty \Psi U_1^n\Psi^*(\sH'\vee\Psi\tau_1\sH)\\
&=& \Psi\bigvee_{n=0}^\infty
U_1^n((\Psi^*\sH')\vee\Psi^*\Psi\tau_1\sH)=
\Psi\bigvee_{n=0}^\infty U_1^n(\sH'\vee\tau_1\sH).
\end{eqnarray*}

We say that an isometric coupling $\coupon{}$ of $\lifset$ is
\emph{special} if $\sK$ is a Hilbert direct sum of the space
$\sH'$, the space $\sD_A$ and some Hilbert space $\sX$, that is,
$\sK=\sH'\oplus\sD_A\oplus \sX$, and the action of $\tau$ is given
by $\tau h= Ah\oplus D_A h\oplus 0$, where $0$ is the zero vector
in $\sX$. In other words, an isometric coupling $\coupon{}$ of
$\lifset$ is special  if, in addition,  $\sD_A$ is a subspace of
$\sM$, where $\sM=\sK\ominus\sH'$, and $\tau$ admits a matrix
representation of the form
\[
\tau=\mat{c}{A\\\Pi_{\sD_A}^*D_A}:\sH\to\mat{c}{\sH'\\\sM}.
\]

The importance of special isometric couplings follows from Theorem
\ref{thm:ic-schur} below. To prove this theorem we need a few
auxiliary propositions. The  first also settles the question of
existence of special isometric couplings.

\begin{proposition}\label{prop:special}
Every isometric coupling is isomorphic to a special isometric
coupling.
\end{proposition} \bpr Let $\coupon{}$ be an isometric
coupling of $\lifset$, and put $\sM=\sK\ominus\sH'$. Since $\tau$
is an isometry and $\Pi_{\sH'}\tau=A$, the operator $\tau$ admits
a matrix representation of the form:
\begin{equation}
\label{notregtau}\tau=\mat{c}{A\\\Gamma D_A}:\sH\to \mat{c}{\sH'\\
\sM}  \mbox{ where } \Gamma:\sD_A\to\sM \mbox{ is an isometry};
\end{equation}
see Section IV.1 of \cite{ff} or Section XXVII.5 of~\cite{ggk}.
Now  let $\sD=\im \Gamma$,  and put $\sX=\sM\ominus \sD$. Then
$\sD$ is closed, and we can view $\Gamma$ as a unitary operator
from $\sD_A$ onto $\sD$. Define the unitary operator $\sigma$ by
\[
\sigma=\mat{ccc}{I_{\sH'}&0&0\\0&\Gamma&0\\0&0&I_\sX}
:\mat{c}{\sH'\\\sD_A\\\sX}\to\mat{c}{\sH'\\\sD\\\sX}.
\]
Also define $U_0=\sigma^*U\sigma$ and $\tau_0=\sigma^*\tau$. Then
$\coup{0}$ is a special isometric coupling of $\{T',A\}$ which is
isomorphic to $\coup{}$. \epr

\medskip Since minimality of  isometric couplings is preserved under
isomorphisms, and isometric couplings do exist (see the third
paragraph of this section), the above proposition shows that any
$\lifset$ admits a special minimal isometric coupling.

Recall that an isometric lifting $U$ of $T'$ can always be
represented (see~(\ref{eq:U1})) in the following form:
\begin{equation}\label{eq:U3}
U = \left[\begin{array}{cc}
  T' & 0 \\
  Y_1D_{T'} & Y_2
\end{array}\right]\mbox{ on }\left[\begin{array}{c}
  \sH' \\
  \sM
\end{array}\right] \mbox{ where \ } Y=\mat{cc}{Y_1&Y_2}: \left[\begin{array}{c}
 \sD_{T'} \\
  \sM
\end{array}\right] \to \sM
\end{equation}
is an isometry. According to  (\ref{eq:U2})  this $U$ also admits
a matrix representation  of the form:
\[
U=\mat{cc}{U'&0\\0&\tilde{U}} \mbox{ on }\mat{c}{\sK'\\
\tilde{\sK}}\quad \mbox{where}\quad\sK'=\bigvee_{n=0}^\infty
U^n\sH',\quad \tilde{\sK}=\sK\ominus\sK'.
\]
Here $U'$ on $\sK'$ is the  minimal isometric lifting of $T'$
associated with $U$ (see Section~\ref{sec:IL}), and $\tilde{U}$ is
an isometry on $\tilde{\sK}$. We can now state the next
proposition.

\begin{proposition}\label{prop:minimal}
Let $\{U,\tau\}$ be an isometric coupling of $\lifset$, where $U$
is determined by  $(\ref{eq:U3})$  and $\tau$ by
$(\ref{notregtau})$. Set $\sD=\im \Gamma$, where $\Gamma$ is given
by $(\ref{notregtau})$, and for $Y$ in  $(\ref{eq:U3})$ consider
the following operator matrix representation:
\begin{equation}\label{matY2} Y =
\left[
\begin{array}{cc}
  D & C \\
  B & Z \\
\end{array}\right]:
\left[\begin{array}{c}
  \sD_{T'}\oplus\sD \\
  \sX
\end{array}\right]\rightarrow
\left[\begin{array}{c}
  \sD \\
  \sX
\end{array}\right]
\quad \mbox{where}\quad \sX=\sM\ominus\sD.
\end{equation}
Then $\{U,\tau\}$ is a minimal  isometric coupling of $\lifset$ if
and only if the pair $\{Z,B\}$ is controllable.
\end{proposition}
\bpr Since $\tau$ is given by (\ref{notregtau}), the space
$\sH'\oplus\tau \sH$ is equal to $\sH'\oplus \sD$. Thus we have to
show that $\sH'\oplus \sD$ is cyclic for $U$ if and only if the
pair $\{Z, B\}$ is controllable. To do this we associate with $U$
two auxiliary operators, namely
\[
\check{U}= \left[\begin{array}{cc}
  0 & 0 \\
  Y_1D_{T'} & Y_2
\end{array}\right]\mbox{ on }\left[\begin{array}{c}
  \sH' \\
  \sM
\end{array}\right], \quad M=\left[\begin{array}{cc}
  0 & 0 \\
  Y_1& Y_2
\end{array}\right]\mbox{ on }\left[\begin{array}{c}
  \sD_{T'} \\
  \sM
\end{array}\right].
\]
Notice that the range of $U-\check{U}$ belongs to $\sH'$. Since
$\sH'\subset \sH'\oplus \sD$, this implies that $\sH'\oplus \sD$
is cyclic for $U$ if and only if $\sH'\oplus \sD$ is cyclic for
$\check{U}$. By induction one proves that for $n=1,2,3, \ldots$
we have
\[
\check{U}^n= \left[\begin{array}{cc}
  0 & 0 \\
  Y_2^{n-1}Y_1D_{T'} & Y_2^n
\end{array}\right]\mbox{ on }\left[\begin{array}{c}
  \sH' \\
  \sM
\end{array}\right], \quad M=\left[\begin{array}{cc}
  0 & 0 \\
  Y_2^{n-1}Y_1& Y_2^n
\end{array}\right]\mbox{ on }\left[\begin{array}{c}
  \sD_{T'} \\
  \sM
\end{array}\right].
\]
On the other hand
\[
\overline{\left[\begin{array}{cc}
  Y_2^{n-1}Y_1D_{T'}& Y_2^n
\end{array}\right]\left[\begin{array}{c}
  \sH' \\
  \sD
\end{array}\right]}=\overline{\left[\begin{array}{cc}
  Y_2^{n-1}Y_1& Y_2^n
\end{array}\right]\left[\begin{array}{c}
  \sD_{T'} \\
  \sD
\end{array}\right]},\quad
 n= 1,2,3, \ldots,
\]
and hence $\sH'\oplus \sD$ is cyclic for $\check{U}$ if and only
if $\sD_{T'}\oplus \sD$ is cyclic for $M$.

It remains to prove that $\sD_{T'}\oplus \sD$ is cyclic for $M$ if
and only if the pair $\{Z,B\}$ is controllable. Using
$Y=\mat{cc}{Y_1&Y_2}$ and (\ref{matY2}), we see that
\[
M=\left[\begin{array}{cc}
  \check{D} & \check{C} \\
  B& Z
\end{array}\right],
\]
where $\check{D}=\Pi_{\sD}^*D$ and $\check{C}=\Pi_{\sD}^*C$ with
$\Pi_{\sD}$ equal to the orthogonal projection of
$\sD_{T'}\oplus\sD$ onto $\sD$. By employing Lemma \ref{lem:cont}
with $\sU = \sY =\sD_{T'}\oplus \sD$, we see that $\sD_{T'}\oplus
\sD$ is cyclic for $M$ if and only if the pair $\{Z, B\}$ is
controllable, which completes the proof.\epr

\begin{proposition}
\label{propisocoupl}  Let  $\coupon{1}$ and $\coupon{2}$ be
special isometric couplings of $\lifset$.  For $j=1,2$ set
$\sX_j=\sK_j\ominus(\sH'\oplus\sD_A)$, and let  $Y(j)$ be the
isometry from $\sD_{T'}\oplus\sD_A\oplus\sX_j$ into
$\sD_A\oplus\sX_j$ corresponding to $U_j$ via
\textup{(\ref{eq:U3})}. Consider the following operator matrix
representation
\[
Y(j)=\mat{cc}{D_j&C_j\\B_j&Z_j}
:\mat{c}{\sD_{T'}\oplus\sD_A\\\sX_j}\to\mat{c}{\sD_A\\\sX_j}
\qquad \mbox{for } j=1,2.
\]
Then $\coup{1}$ and $\coup{2}$ are isomorphic if and only if
$\{Z_1,B_1,C_1,D_1\}$ and $\{Z_2,B_2,C_2,D_2\}$ are unitarily
equivalent realizations.
\end{proposition}
\bpr Assume that $\{Z_1,B_1,C_1,D_1\}$ and $\{Z_2,B_2,C_2,D_2\}$
are unitarily equivalent, that is, $D_1=D_2$ and there exists a
unitary operator  $W$ from $\sX_1$ onto $\sX_2$ such that
\begin{equation}\label{eq:phix1}
W Z_1=Z_2W,\quad W B_1=B_2\quad\mbox{and}\quad C_1=C_2W.
\end{equation}
Now let $\Phi$ be the unitary operator from $\sK_1$ onto $\sK_2$
defined by
\begin{equation}\label{eq:phix}
\Phi=\mat{ccc}{I_{\sH'}&0&0\\
0&I_{\sD_A}&0\\0&0&W} :\mat{c}{\sH'\\
\sD_A\\ \sX_1}\to\mat{c}{\sH'\\\sD_A\\\sX_2}.
\end{equation}
Then $\Phi h=h$ for all $h$  in $\sH'$.  Because  $\coup{1}$ and
$\coup{2}$ are special, we see  that
\begin{equation}\label{eq:taui}
\tau_j=\mat{c}{A\\ D_A\\0}:\sH\to\mat{c}{\sH'\\\sD_A\\\sX_j}
\qquad \mbox{for }j=1,2.
\end{equation}
Hence $\Phi\tau_1=\tau_2$. Using the appropriate operator matrix
decomposition we arrive at
\begin{eqnarray}
\Phi U_1 &=&\mat{ccc}{I_{\sH'}&0&0\\0&I_{\sD_A}&0\\0&0&W}
\mat{ccc}{T'&0&0\\
D_1\Pi_{\sD_{T'}}^*D_{T'}&D_1\Pi_{\sD_{A}}^*&C_1\\
B_1\Pi_{\sD_{T'}}^*D_{T'}&B_1\Pi_{\sD_{A}}^*&Z_1} \nonumber\\
&=&\mat{ccc}{T'&0&0\\
D_1\Pi_{\sD_{T'}}^*D_{T'}&D_1\Pi_{\sD_{A}}^*&C_1\\
WB_1\Pi_{\sD_{T'}}^*D_{T'}&WB_1\Pi_{\sD_{A}}^*&WZ_1}.
\label{eq:calc1}
\end{eqnarray}
A similar calculation shows that
\begin{eqnarray}
U_2\Phi
&=&\mat{ccc}{T'&0&0\\
D_2\Pi_{\sD_{T'}}^*D_{T'}&D_2\Pi_{\sD_{A}}^*&C_2\\
B_2\Pi_{\sD_{T'}}^*D_{T'}&B_2\Pi_{\sD_{A}}^*&Z_2}
\mat{ccc}{I_{\sH'}&0&0\\0&I_{\sD_A}&0\\0&0&W}\nonumber\\
&=&\mat{ccc}{T'&0&0\\
D_2\Pi_{\sD_{T'}}^*D_{T'}&D_2\Pi_{\sD_{A}}^*&C_2W\\
B_2\Pi_{\sD_{T'}}^*D_{T'}&B_2\Pi_{\sD_{A}}^*&Z_2W}.\label{eq:calc2}
\end{eqnarray}
Because $D_1=D_2$ and  (\ref{eq:phix1}) holds, we see that $\Phi
U_1=U_2\Phi$. In other words,  $\coup{1}$ and $\coup{2}$ are
isomorphic.

Conversely assume that $\coup{1}$ and $\coup{2}$ are isomorphic.
Then there exists a  unitary operator $\Phi$ from $\sK_1$ onto
$\sK_2$ such that $\Phi h=h$ for all $h$  in $\sH'$ and
$\Phi\tau_1=\tau_2$ and $\Phi U_1=U_2\Phi$. Because   $\tau_1$ and
$\tau_2$ admit matrix representations of the form presented in
(\ref{eq:taui}) and $\im D_A$ is dense in $\sD_A$, we see that
$\Phi h=h$ for all $h\in\sD_A$. So $\Phi$ admits a matrix
representation as in (\ref{eq:phix})  where $W$ is a unitary
operator from $\sX_1$ onto $\sX_2$.   By combining   $\Phi
U_1=U_2\Phi$ with the matrix representations for  $\Phi U_1$ in
(\ref{eq:calc1}) and $U_2\Phi$ in  (\ref{eq:calc2}), we see   that
\[
D_1=D_2,\quad W Z_1=Z_2W,\quad W B_1=B_2\quad\mbox{and}\quad
C_1=C_2W.
\]
Hence $\{Z_1,B_1,C_1,D_1\}$ and $\{Z_2,B_2,C_2,D_2\}$ are
unitarily equivalent realizations. \epr

\begin{theorem}
\label{thm:ic-schur} Let $\lifset$ be a pair of contractions, $T'$
acting on $\sH'$ and $A$ from $\sH$ into $\sH'$. Then there is a
one to one map from the set of minimal isometric couplings of
$\lifset$, with isomorphic ones being identified, onto the Schur
class $\eS(\sD_A,\sD_{T'}\oplus\sD_A)$. This map is defined as
follows. Let $\{U,\tau\}$ be a minimal isometric coupling  of
$\lifset$, which   may be assumed to be special, by Proposition
\textup{\ref{prop:special}}. Define
\begin{equation}\label{eq:ICF}
F_{\{U,\tau\}}(\la)= \Pi_{\sD_{T'}\oplus\sD_A}Y^*(I_{\sM}-\la
J_\sX^\prime Y^*)^{-1}\Pi^*_{\sD_A},
\end{equation}
where $Y$ is the isometry uniquely determined by $U$ via
$(\ref{eq:U3})$, $\sX=\sM\ominus \sD_A$, and $J_\sX^\prime$ is the
partial isometry from $\sD_{T'}\oplus\sD_A\oplus\sX$ to
$\sD_A\oplus\sX$ given by
\[
J_\sX^\prime=\mat{cc}{0&0\\0&I_{\sX}}:
\mat{c}{\sD_{T'}\oplus\sD_A\\\sX}\to\mat{c}{\sD_A\\\sX}.
\]
Then $\{U,\tau\} \mapsto F_{\{U,\tau\}}$ is the desired map.
\end{theorem}

\bpr We know from Proposition~\ref{prop:special} that every
isometric coupling is isomorphic to a special one. So without loss
of generality we can assume the isometric couplings to be special.

{}From Proposition~\ref{prop:minimal} and Section~\ref{sec:IL} it
is clear that there is a one to one correspondence between the
special minimal isometric couplings of $\lifset$ and the
isometries $Y$ mapping the space $\sD_{T'}\oplus\sD_A\oplus\sX$
into $\sD_A\oplus\sX$,  where $\sX$ is some Hilbert space and
 the pair $\{\Pi_\sX Y\Pi_\sX^*,\Pi_\sX Y\Pi_{\sD_{T'}\oplus\sD_A}^*\}$ is
controllable. In fact, this one to one correspondence is provided
by (\ref{eq:U3}). Furthermore, formula (\ref{matY2}) establishes a
one to one correspondence between the isometries $Y$ mapping the
space $\sD_{T'}\oplus\sD_A\oplus\sX$ into $\sD_A\oplus\sX$ and the
isometric realizations
\[ \{Z,B,C,D;\sX,\sD_{T'}\oplus\sD_A, \sD_A\},
\]
and in this one to one correspondence $Z=\Pi_\sX Y\Pi_\sX^*$ and
$B=\Pi_\sX Y\Pi_{\sD_{T'}\oplus\sD_A}^*$. {}From
Theorem~\ref{th:CAF-IR} we know that there is a one to one
correspondence between   the controllable isometric realizations,
with the unitarily equivalent ones being identified, and the
$\eS(\sD_{T'}\oplus\sD_A, \sD_A)$  Schur class functions. Next,
note that the map $G\mapsto F$, where
$F(\la)=G(\overline{\la})^*$, is a one to one map from
$\eS(\sD_{T'}\oplus\sD_A, \sD_A)$ onto $\eS(\sD_A,
\sD_{T'}\oplus\sD_A)$. Following up all these one to one
correspondences and using the results of Section \ref{sec:IR} we
see that the map from a special minimal isometric coupling
$\{U,\tau\}$ to $F$ is given by $F=F_{\{U,\tau\}}$. To complete
the proof, it remains to apply Proposition \ref{propisocoupl}.
\epr

\medskip We conclude this section with a lemma that will be useful
in the next section.

\begin{lemma}\label{lemisocoup12}
Let $\coupon{1}$ and $\coupon{2}$ be isomorphic isometric
couplings of $\lifset$,  and let $V$ on $\sH'\oplus H^2(\sD_{T'})$
be the Sz.-Nagy-Sch\"affer minimal isometric lifting of $T'$. For
$j=1,2$ let $\Phi_j$ be the unique isometry associated with $T'$
intertwining $V$ and $U_j$. Then
\[
\Phi_1^*\tau_1=\Phi_2^*\tau_2.
\]
\end{lemma}
\bpr Let $\Psi$ from $\sK_1$ to $\sK_2$ be an isomorphism from
$\coup{1}$ to $\coup{2}$. Define $\Theta$ from $\sH'\oplus
H^2(\sD_{T'})$ into $\sK_1$ by setting $\Theta=\Psi^*\Phi_2$. Then
$\Theta$ is an isometry, $\Theta h=\Psi^*\Phi_2 h=\Psi^* h=h$ for
all $h\in\sH'$, and
\[
\Theta V=\Psi^*\Phi_2 V=\Psi^*
U_2\Phi_2=U_1\Psi^*\Phi_2=U_1\Theta.
\]
So, by Theorem~\ref{thmlift} (see also the last paragraph of
Section \ref{sec:IL}), the operator $\Theta$ is the unique
isometry associated with $T'$ intertwining $V$ and $U_1$, that is,
$\Theta=\Phi_1$. It follows that $\Phi_1=\Psi^*\Phi_2$, and hence
$ \Phi_1^*\tau_1=\Phi_2^*\Psi\tau_1=\Phi_2^*\tau_2$, which
completes the proof.\epr

\section{Main theorem for the case when
$R^*R=Q^*Q$}\label{sec:RCLT0} \setcounter{equation}{0}

In this section $\{A,T',V, R,Q\}$ is a lifting data set, with $V$
on $\sH'\oplus H^2(\sD_{T'})$ being the Sz.-Nagy-Sch\"affer
minimal isometric lifting of $T'$. In particular, $T'AR=AQ$ and
$R^*R\leq Q^*Q$. Recall that  $B$ from $\sH$ into $\sH'\oplus
H^2(\sD_{T'})$ is a  contractive interpolant  for $\liftsetV$ if
$B$ is a contraction satisfying $\Pi_{\sH'}B=A$ and $VBR=BQ$.

Our aim is to prove Theorem \ref{mainth} assuming that
$R^*R=Q^*Q$. First let us reformulate Theorem \ref{mainth} for
this case. For this purpose note that for  $R^*R=Q^*Q$ the  spaces
$\sF$ and $\sF'$ defined by (\ref{spacesF}) are given by
\[
\sF=\overline{D_AQ\sH_0}
 \quad \mbox{ and } \quad
\sF'=\overline{\left[
\begin{array}{c}
D_{T'}AR\\D_AR
\end{array}
\right]\sH_0}.
\]
Observe  that $\sF \subset \sD_A$ and $\sF' \subset \sD_{T'}\oplus
\sD_A$. Furthermore, the unitary operator $\omega$ mapping $\sF$
onto $\sF'$ in (\ref{omega}) is now determined by
\begin{equation}
\label{omega2} \omega(D_AQh)=\left[
\begin{array}{c}
D_{T'}AR\\D_AR
\end{array}
\right]h, \quad h\in \sH_0.
\end{equation}
The following is the main result of this section.

\begin{theorem}
\label{mainth2} Let $\liftsetV$ be a lifting data set, where $V$
on $\sH'\oplus H^2(\sD_{T'})$ is the Sz.-Nagy-Sch\"affer minimal
isometric lifting of $T'$, and assume that $R^*R = Q^*Q$. Then all
contractive interpolants  $B$ for $\liftsetV$ are given by
\begin{equation}
\label{sols2}Bh=\left[
\begin{array}{c}
Ah\\ \Pi_{T'}F(\lambda)(I-\lambda \Pi_A F(\lambda))^{-1}D_Ah
\end{array}
\right], \quad h\in \sH,
\end{equation}
where $F$ is any  function  in $\eS(\sD_A, \sD_{T'}\oplus \sD_A)$
satisfying $F(0)|\sF=\omega$. Here $\omega$ is the unitary
operator defined in $(\ref{omega2})$ while  $\Pi_{T'}$ and $\Pi_A$
are the projections given by
\[\Pi_{T'}=\left[\begin{array}{cc} I&0
\end{array}\right]: \left[
\begin{array}{c}
\sD_{T'}\\ \sD_A
\end{array}
\right]\to \sD_{T'}
\quad \mbox{and} \quad\Pi_A=\left[\begin{array}{cc} 0&I
\end{array}\right]: \left[
\begin{array}{c}
\sD_{T'}\\ \sD_A
\end{array}
\right]\to \sD_A.
\]
\end{theorem}

The proof of the above theorem will be based on a further
refinement (which we present in two propositions) of the theory of
isometric couplings presented in the previous section. In fact, to
obtain contractive interpolants for $\liftsetV$ we shall need
isometric couplings $\{U,\tau\}$ of $\{T',A\}$ satisfying the
additional intertwining relation $U\tau R=\tau Q$. This is the
contents of the first proposition (Proposition \ref{prop:couplsol}
below). The existence of such couplings is guaranteed by the
second proposition (Proposition \ref{ic-schurRQ} below), which is
based on Theorem \ref{thm:ic-schur}. In the sequel, for
simplicity, we shall write $\sV$ for the space $\sH'\oplus
H^2(\sD_{T'})$.

\begin{proposition}
\label{prop:couplsol} Let $\liftsetV$ be a lifting data set, with
$V$ on $\sV$ being the Sz.-Nagy-Sch\"affer minimal isometric
lifting of $T'$, and assume that $R^*R = Q^*Q$. Let $\{U\mbox{ on
}\sK,\tau\}$ be an isometric coupling of $\{T',A\}$ satisfying
$U\tau R=\tau Q$, and let $\Phi$ be the unique isometry from $\sV$
into $\sK$ associated with $T'$ intertwining $V$ with $U$. Then
\begin{equation}
\label{couplsol} B=\Phi^*\tau
\end{equation}
is a contractive interpolant for $\liftsetV$, and all contractive
interpolants for this data set are obtained in this way. More
precisely, if $B$ a contractive interpolant for $\liftsetV$, then
there exists a minimal special  isometric coupling $\{U,\tau\}$ of
$\{T',A\}$ such that $B=\Phi^*\tau$ and $U\tau R=\tau Q$.
\end{proposition}
\bpr First let us show that $B$ defined by (\ref{couplsol}) is a
contractive interpolant for the data set $\{A,T',V,R,Q\}$.
Obviously, $B$ is a contraction. Put  $\sK'=\im \Phi$. Recall (see
Section~\ref{sec:IL}) that $\Phi\Phi^*$ is the orthogonal
projection of $\sK$ onto $\sK'$. {}From Theorem \ref{thmlift} we
know that $\sK'=\bigvee_{n\geq0}U^n\sH'$ is a reducing subspace
for $U$. It follows that $U$ commutes with $\Phi\Phi^*$. Since
$\Phi^*\Phi$ is the identity operator on $\sV=\sH'\oplus
H^2(\sD_{T'})$, we obtain
\[
VBR=\Phi^*U\Phi(\Phi^*\tau)R=\Phi^*(U\Phi\Phi^*)\tau
R=\Phi^*\Phi\Phi^*U\tau R=\Phi^*\tau Q=BQ.
\]
Thus $B$ is a contractive interpolant for  $\{A,T',V,R,Q\}$.

To prove the reverse implication, assume that $B$ is a contractive
interpolant. We have to construct a minimal special isometric
coupling $\{U,\tau\}$ of $\{T',A\}$ satisfying $U\tau R=\tau Q$
such that $B$ is given by (\ref{couplsol}). Since $B$ is a
contraction, we may consider the  subspaces
\[
\tilde\sF=\overline{D_BR\sH_0} \quad \mbox{and}\quad
\tilde\sF'=\overline{D_BQ\sH_0}.
\]
Using $VBR=BQ$ with $R^*R=Q^*Q$, and the fact that $V$ is an
isometry, we see that for each $h\in\sH_0$ we have
\begin{eqnarray*}
\|D_BQh\|^2 &=&\|Qh\|^2-\|BQh\|^2 =\|Rh\|^2-\|VBRh\|^2\\
&=&\|Rh\|^2-\|BRh\|^2 =\|D_BRh\|^2.
\end{eqnarray*}
Hence there exists a unique unitary operator $\tilde\oo$ from
$\tilde \sF$ onto $\tilde\sF'$ such that $\tilde\oo D_BR=D_BQ$.
Next, define
 the subspaces
\[
\tilde \sG=\sD_B\ominus\tilde \sF \quad \mbox{and}\quad
\tilde\sG'=\sD_B\ominus\tilde \sF'.
\]
Notice  that $\sD_B=\tilde \sF\oplus \tilde \sG$ and $\sD_B=\tilde
\sF' \oplus \tilde \sG'$.  Thus $\tilde\oo$ defines a partial
isometry $\om$ on $\sD_B$ as follows:
\[
\om=\left[\begin{array} {cc} \tilde\oo&0\\0&0
\end{array}\right]:\left[\begin{array} {c}
\tilde\sF\\\tilde\sG
\end{array}\right]\to\left[\begin{array} {c}
\tilde\sF'\\\tilde\sG'
\end{array}\right].
\]
Observe that $\sD_\Omega$ coincides with $\tilde\sG$. Define
$V_\om$ to be the Sz.-Nagy-Sch\"affer minimal isometric lifting of
$\om$ on $\sV_\om=\sD_B\oplus H^2(\tilde\sG)$. Thus $V_\om$ has
the following  operator matrix representation
\[
V_\om=\left[\begin{array} {ccc} \tilde\oo&0&0\\0&0&0\\
0&E_{\tilde\sG}&S_{\tilde\sG}
\end{array}\right]:\left[\begin{array} {c}
\tilde\sF\\\tilde\sG\\H^2(\tilde\sG)
\end{array}\right]\to\left[\begin{array} {c}
\tilde\sF'\\\tilde\sG'\\H^2(\tilde\sG)
\end{array}\right].
\]
Here $E_{\tilde\sG}$ is the canonical embedding of $\tilde\sG$
onto the space of constant functions in $H^2(\tilde\sG)$, and
$S_{\tilde\sG}$ is the unilateral shift on the Hardy space
$H^2(\tilde\sG)$. Since $V_\om$ is a minimal isometric lifting of
$\om$, we have $\sV_\om=\sD_B\oplus
H^2(\tilde\sG)=\bigvee_{n=0}^\infty V_\om^n\sD_B$. Now, put
\begin{eqnarray}
&&U_\om=\left[\begin{array} {cc} V&0\\0&V_\om
\end{array}\right]\ \textup{on}\ \left[\begin{array} {c}
\sV\\ \sV_\om
\end{array}\right],\label{Uomega} \\
\noalign{\vskip6pt} &&\tau_\om=\left[\begin{array} {c}
B\\\Pi_{\sD_B}^*D_B
\end{array}\right]:\sH\to\left[\begin{array} {c}
\sV\\ \sV_\om
\end{array}\right].\label{tauomega}
\end{eqnarray}
Since $V$ on $\sV=\sH'\oplus H^2(\sD_{T'})$ is an isometric
lifting of $T'$, the operator $U_\om$ is an isometric lifting of
$T'$, and $A=\Pi_{\sH'}\tau_\Omega$. It follows that $\{U_\om,
\tau_\om\}$ is an isometric coupling of $\lifset$. Notice that
$\sV\vee\tau_\om\sH=\sV\oplus \sD_B$. Because $\sH'$ is cyclic for
$V$ and $\sD_B$ is cyclic for $V_\Omega$, the reducing
decomposition of $U_\Omega$ in (\ref{Uomega}) shows that
$\sH'\vee\tau_\om\sH$ is cyclic for $U_\Omega$. In other words,
the isometric coupling $\{U_\om,\tau_\om\}$ is minimal. Since
$VBR=BQ$, the construction of $U_\om$ and $\tau_\om$ implies that
\begin{equation}\label{intertwURQ}
U_\om\tau_\om R=\tau_\om Q \quad \mbox{and}\quad
B=\Pi_{\sV}\tau_\om.
\end{equation}
Indeed, for $h\in\sH_0$ we have
\[
U_\om\tau_\om Rh=\left[\begin{array} {cc} V&0\\0&V_\om
\end{array}\right]\left[\begin{array} {c}
BRh\\D_BRh
\end{array}\right]=\left[\begin{array} {c}
VBRh\\ V_\om D_BRh
\end{array}\right].
\]
However, $D_BR\sH_0\subset\tilde\sF$, and hence $V_\om
D_BRh=\tilde\oo D_BRh=D_BQh$, which follows from  the definition
of $\tilde\oo$. Since, by assumption, $VBRh=BQh$, we see that
\[
U_\om\tau_\om Rh=\left[\begin{array} {c} BQh\\D_BQh
\end{array}\right]=\tau_\om Qh,
\]
which proves the first identity in (\ref{intertwURQ}). The second
is clear from the definition of $\tau_\om$.

{}From the construction of $U_\om$ it follows that the unique
isometry $\Phi_\om$ associated with $T'$ that intertwines $V$ with
$U_\om$ is equal to $\Pi_\sV^*$, where $\Pi_\sV$ is the orthogonal
projection of $\sV\oplus \sV_\om$ onto $\sV$. This together with
the second identity in (\ref{intertwURQ}) yields
\begin{equation}
\label{couplsolo} B=\Pi_\sV\tau_\om=\Phi_\om^*\tau_\om.
\end{equation}

By Proposition \ref{prop:special} and the fact that minimality of
isometric couplings is preserved under isomorphisms, there exists
a minimal special isometric coupling $\{U,\tau\}$ of $\{T',A\}$
which is isomorphic to $\{U_\om,\tau_\om\}$. Using Lemma
\ref{lemisocoup12} and formula (\ref{couplsolo}) we obtain
$B=\Phi^* \tau$, where $\Phi$ is the unique isometry associated
with $T'$ that intertwines $V$ with $U$.

It remains to prove that $U\tau R=\tau Q$. Let $\Psi$ be the
isomorphism that transforms $\{U_\om,\tau_\om\}$ into
$\{U,\tau\}$. In particular, $\Psi \tau_\om= \tau$. Moreover
formula (\ref{eq:morph3a}) yields $\Psi U_\om\tau_\om=U\tau$.
Since $U_\om\tau_\om R=\tau_\om Q$, it follows that $U\tau R=\Psi
U_\om\tau_\om R=\Psi \tau_\om Q= \tau Q$. \epr

\begin{proposition}
\label{ic-schurRQ} Let $\liftsetV$ be a lifting data set. Assume
that $R^*R = Q^*Q$, and let $\omega$ be the unitary operator
defined in
 $(\ref{omega2})$.
Consider a minimal special isometric coupling $\{U,\tau\}$ of
$\lifset$, and let $F_{\{U,\tau\}}$
 be the function in the Schur class $\eS(\sD_A,\sD_{T'}\oplus\sD_A)$
 defined by $(\ref{eq:ICF})$. Then $U\tau R=\tau Q$ if and only
$F_{\{U,\tau\}}(0)|\sF=\omega$.  In particular, there exists a
special isometric coupling $\{U,\tau\}$ of $\{T',A\}$ satisfying
$U\tau R=\tau Q$.
\end{proposition}

It will be convenient first to prove  the following lemma.

\begin{lemma}
\label{lem:RQconstr} Let $\liftset$ be a lifting data set
satisfying $R^*R = Q^*Q$.  Let $\{U,\tau\}$ be a special isometric
coupling of $\{T',A\}$, and consider its operator matrix
representation of the form
\begin{equation}\label{eq:U3a}
U = \left[\begin{array}{cc}
  T' & 0 \\
  Y_1D_{T'} & Y_2
\end{array}\right]\mbox{ on }\left[\begin{array}{c}
  \sH' \\
  \sM
\end{array}\right]\mbox{ where } Y=\mat{cc}{Y_1&Y_2}: \left[\begin{array}{c}
 \sD_{T'} \\
  \sM
\end{array}\right] \to \sM
\end{equation}
is an isometry. Then $U\tau R=\tau Q$ if and only if
$Y|\sF'=\omega^*$.
\end{lemma}
\bpr Since the coupling is special, the space $\sD_A$ is a
subspace of $\sM$ and
\[
\tau=\mat{c}{A\\ \Pi_{\sD_A}^*D_A}:\sH\to\mat{c}{\sH'\\\sM}.
\]
It follows that for $h$ in $\sH_0$, we have
\begin{eqnarray*}
U\tau Rh &=&\mat{cc}{T'&0\\
Y_1D_{T'}&Y_2}\mat{c}{ARh\\  D_ARh}= \mat{c}{T'A
Rh\\Y_1D_{T'}ARh+Y_2D_{A}Rh}\\
&=&\mat{c}{AQh\\Y\mat{c}{D_{T'}ARh\\D_ARh}}.
\end{eqnarray*}
Thus
\[
U\tau R=\tau Q \quad \Longleftrightarrow\quad
Y\mat{c}{D_{T'}AR\\D_AR}h=D_AQh,  \quad h\in \sH_0.
\]
Since $Y$ is an isometry, we see that $U\tau R=\tau Q$ if and only
if $Y|\sF'$ is a unitary operator from $\sF'$ onto $\sF$ with the
same action as $\omega^*$. Because of the uniqueness of $\omega$,
this proves the lemma. \epr

\medskip\noindent\textbf{Proof of Proposition \ref{ic-schurRQ}.}
Let $Y$ be the isometry determined by the operator matrix
representation for $U$ in (\ref{eq:U3a}), and set
$F=F_{\{U,\tau\}}$. {}From Lemma \ref{lem:RQconstr} we know that
$U\tau R=\tau Q$ if and only if $Y|\sF'=\oo^*$. Thus we have to
show that
\begin{equation}
\label{equivcondFY} F(0)|\sF=\oo\Longleftrightarrow Y|\sF'=\oo^*.
\end{equation}
By consulting  (\ref{eq:ICF}) we see that
$F(0)=\Pi_{\sD_{T'}\oplus\sD_A}Y^*\Pi_{\sD_A}^*$. The fact that
$\sF'\subset \sD_{T'}\oplus \sD_A$ allows us to view $\oo$ as an
isometry from $\sF$ into $\sD_{T'}\oplus \sD_A$. Since
$\sF\subset\sD_A$, it follows that the first condition in
(\ref{equivcondFY}) is equivalent to
\[
Y^*|\sF=\left[\begin{array} {c} \oo\\\gamma
\end{array}\right]:\sF\to\left[\begin{array} {c}
\sD_{T'}\oplus\sD_A\\ \sX
\end{array}\right],
\]
where $\gamma$ is some operator from $\sF$ into $\sX$.  However,
$\oo$ is an isometry and $Y^*|\sF$ is a contraction. This implies
that $\gamma=0$. We conclude that the first condition in
(\ref{equivcondFY}) is equivalent to $Y^*|\sF=\oo$. By taking
adjoints, and using that $\oo$ is a unitary operator from $\sF$
onto $\sF'$, we see that the same holds true for the second
condition in (\ref{equivcondFY}).

Finally, for each $\la \in \BD$ define the operator $F(\la)$ from
$\sD_A$ into $\sD_{T'}\oplus \sD_A$ by setting $F(\la)d=\oo\Pi_\sF
d$ for $d\in \sD_A$. Here $\Pi_\sF$ is the orthogonal projection
of $\sD_A$ onto $\sF$. Then $F$ belongs to the Schur class
$\eS(\sD_A,\sD_{T'}\oplus\sD_A)$. Hence, by Theorem
\ref{thm:ic-schur}, there exists a minimal special isometric
coupling $\{U,\tau\}$ of $\lifset$ such $F=F_{\{U,\tau\}}$. Since
$F(0)|\sF=\oo$, we conclude that $U\tau R=\tau Q$.\epr

\medskip\noindent\textbf{Proof of Theorem \ref{mainth2}.} We split
the proof into two parts.

\smallskip\noindent\textbf{Part 1.} Assume   that  $B$ is a contractive
interpolant for the data set $\{A,T',V,R,Q\}$. Since $R^*R=Q^*Q$,
we know from Proposition \ref{prop:couplsol} that there exists a
(minimal) special isometric coupling $\{U\ \textup{on}\
\sK,\tau\}$ of $\{T',A\}$ such that $U\tau R=\tau Q$ and
$B=\Phi^*\tau$, where $\Phi$ is the unique isometry from
$\sV=\sH'\oplus H^2(\sD_{T'})$ into $\sK$ associated with $T'$
intertwining $V$ with $U$. Now write
\[
U=\left[\begin{array} {cc} T'&0\\ Y_1D_{T'}&Y_2
\end{array}\right]\textup{ on } \left[\begin{array} {c} \sH'\\
\sM\end{array}\right],
\]
where $\sM=\sK\ominus\sH'$. Since $\{U,\tau\}$ is special, we have
$\sD_A\subset\sM$, and
\[
\tau=\left[\begin{array} {c} A\\ \Pi_{\sD_A}^*D_A
\end{array}\right]:\sH\to\left[\begin{array} {c} \sH'\\
\sM\end{array}\right].
\]
The identity $B=\Phi^*\tau$ and the formula for $\Phi$ in Theorem
\ref{thmlift} show that $B=A\oplus\Lambda^*\Pi_{\sD_A}^*D_A$,
where
\[
(\Lambda^*m)(\la)=Y_1^*(I-\la Y_2^*)^{-1}, \quad m\in\sM\quad
\mbox{and}\quad\la\in\BD.
\]
It follows that
\[
Bh=\left[\begin{array} {c} Ah\\ Y_1^*(I-\la
Y_2^*)^{-1}\Pi_{\sD_A}^*D_Ah \end{array}\right],\quad h\in\sH.
\]
To obtain the expression for $B$ given in (\ref{sols2}) we apply
Proposition \ref{propFY} with $\sD=\sD_A$ and $\sD'=\sD_{T'}$. It
follows that (\ref{sols2}) holds with
$F\in\eS(\sD_A,\sD_{T'}\oplus\sD_A)$ given by
\[ F(\la)=\Pi_{\sD_{T'}\oplus\sD_A}Y^*(I_\sM-\la
J_\sX'Y^*)^{-1},\quad \la\in\BD.
\]
Here $\sX=\sM\ominus\sD_A$, the operator
$\Pi_{\sD_{T'}\oplus\sD_A}$ is the orthogonal projection of
$\sD_{T'}\oplus\sM$ onto $\sD_{T'}\oplus\sD_A$, and
\[
J_\sX'=\left[\begin{array} {cc} 0& P_\sX
\end{array}\right]: \left[\begin{array} {c} \sD_{T'}\\ \sM
\end{array}\right]\to\sM,\quad  Y=\left[\begin{array} {cc} Y_1& Y_2
\end{array}\right]:\left[\begin{array} {c} \sD_{T'}\\ \sM
\end{array}\right]\to\sM.
\]
In other words, using the terminology introduced in Theorem
\ref{thm:ic-schur}, we have $F=F_{\{U,\tau\}}$. Since $U\tau
R=\tau Q$, Proposition \ref{ic-schurRQ} shows that $F(0)|\sF=\oo$,
which completes the first part of the proof.

\smallskip\noindent\textbf{Part 2.} Let
$F$ be any function in $\eS(\sD_A,\sD_{T'}\oplus\sD_A)$ satisfying
$F(0)|\sF=\oo$. We have to show that $B$ defined by (\ref{sols2})
is a contractive interpolant for the given data set. According to
Theorem \ref{thm:ic-schur} there is a minimal special isometric
coupling $\{U,\tau\}$ of $\{T',A\}$ such that $F=F_{\{U,\tau\}}$,
where $F_{\{U,\tau\}}$ is defined by (\ref{eq:ICF}). The fact that
$F(0)|\sF=\oo$ yields $U\tau R=\tau Q$, by Proposition
\ref{ic-schurRQ}.

Since $B$ is given by (\ref{sols2}), we can use Proposition
\ref{propFY} (with $\sD=\sD_A$ and $\sD'=\sD_{T'})$ and Theorem
\ref{thmlift} to show that $B=\Phi^*\tau$, where $\Phi$ is the
unique isometry associated with $T'$ intertwining $V$ (the
Sz.-Nagy-Sch\"affer minimal isometric lifting of $T'$) with $U$.
This allows us to apply Proposition \ref{prop:couplsol} to show
that $B$ is a contractive interpolant. \epr

\section{Proof of the first main theorem} \label{sec:RCLT1}
\setcounter{equation}{0} In this section we shall prove Theorem
\ref{mainth}. The proof will be based on the analogous result for
the case when $R^*R=Q^*Q$, which was proved in the preceding
section, and on Proposition~\ref{propbbtil} below, which allows us
to reduce the general case to the case when $R^*R=Q^*Q$.

Throughout this section $\{A,T',V, R,Q\}$ is a lifting data set
with $V$ being the Sz.-Nagy-Sch\"affer minimal isometric lifting
of $T'$. As before, put $\sD_\circ=\overline{D_\circ\sH_0}$, where
$D_\circ$ is the positive square root of $Q^*Q-R^*R$. Introduce
the following operators:
\begin{eqnarray*}
A_\circ&=&\left[\begin{array} {cc} A&0\\0&I_{\sD_\circ}
\end{array}\right]: \left[\begin{array} {c} \sH\\\sD_{\circ}
\end{array}\right]\to\left[\begin{array} {c} \sH'\\\sD_{\circ}
\end{array}\right],\quad T_\circ'=\left[\begin{array} {cc} T'&0\\0&0
\end{array}\right]\ \textup{on}\ \left[\begin{array} {c} \sH'\\ \sD_{\circ}
\end{array}\right],\\
\noalign{\vskip6pt} R_\circ&=&\left[\begin{array} {c} R\\ D_\circ
\end{array}\right]: \sH_0 \to\left[\begin{array} {c} \sH\\\sD_{\circ}
\end{array}\right],\quad
Q_\circ=\left[\begin{array} {c} Q\\0
\end{array}\right]: \sH_0 \to\left[\begin{array} {c} \sH\\\sD_{\circ}
\end{array}\right],\\
\noalign{\vskip6pt}
\tilde V&=&\left[\begin{array}{cccc} T'&0&0&0\\ 0&0&0&0\\
E_{\sD_{T'}}D_{T'}&0&S_{\sD_{T'}}&0\\
0&E_{\sD_\circ}&0&S_{\sD_\circ}
\end{array}\right]\ \textup{on}\
\left[\begin{array} {c} \sH'\\\sD_{\circ}\\ H^2(\sD_{T'})\\
H^2(\sD_\circ)
\end{array}\right].
\end{eqnarray*}
Here $E_{\sD_{T'}}$ and $E_{\sD_\circ}$ are the canonical
embeddings of $\sD_{T'}$ and $\sD_\circ$ onto the spaces of
constant functions of $H^2(\sD_{T'})$ and $H^2(\sD_\circ)$,
respectively, and $S_{\sD_{T'}}$ and $S_{\sD_A}$ are the forward
shifts on $H^2(\sD_{T'})$ and $H^2(\sD_A)$, respectively.
Identifying $H^2(\sD_{T'}\oplus \sD_\circ)$ with
$H^2(\sD_{T'})\oplus H^2(\sD_\circ)$ it is straightforward to
check that $\tilde{V}$ is the Sz.-Nagy-Sch\"affer minimal
isometric lifting of $T_\circ'$, and that the quintet
\begin{equation}
\label{datatil1} \{A_\circ,T_\circ',\tilde V,R_\circ,Q_\circ\}
\end{equation}
is a lifting data set satisfying
$R_\circ^*R_\circ=Q_\circ^*Q_\circ$.

\begin{proposition}
\label{propbbtil}
 If $\tilde B$ from $\sH\oplus\sD_\circ$ to
$\sH'\oplus\sD_\circ\oplus H^2(\sD_{T'})\oplus H^2(\sD_\circ)$ is
a contractive interpolant for the data set
\textup{(\ref{datatil1})}, then the operator $B$ from $\sH$ to
$\sH'\oplus H^2(\sD_{T'})$, defined by
\begin{equation}
\label{BBtil} B=\Pi_{\sH'\oplus H^2(\sD_{T'})}\tilde B\Pi_\sH^*,
\end{equation}
is a contractive interpolant for the data set $\{A,T,V,R,Q\}$, and
all contractive interpolants for $\{A,T,V,R,Q\}$ are obtained in
this way.
\end{proposition}
\bpr Let $\tilde B$ be a contractive interpolant for the data set
(\ref{datatil1}). Then $\tilde B$ is of the following form
\begin{equation}
\label{Btil1} \tilde B=\left[\begin{array} {cc} A&0\\
0&I_{\sD_\circ}\\ \Gamma_1D_A&0\\ \Gamma_2D_A&0
\end{array}\right] :
\left[\begin{array} {c} \sH\\\sD_{\circ}
\end{array}\right]\to
\left[\begin{array} {c} \sH'\\\sD_{\circ}\\ H^2(\sD_{T'})\\
H^2(\sD_\circ)
\end{array}\right],
\end{equation}
where
\[\left[\begin{array} {c} \Gamma_1\\\Gamma_2
\end{array}\right]:\sD_A\to\left[\begin{array} {c} H^2(\sD_{T'})\\
H^2(\sD_\circ)
\end{array}\right] \textup{is\ a\ contraction}.
\]
Moreover, $\tilde V\tilde BR_\circ=\tilde BQ_\circ$. Now, using
this $\tilde B$, let $B$ be the operator defined by (\ref{BBtil}).
In other words
\[
B=\left[\begin{array} {c} A\\
\Gamma_1 D_A
\end{array}\right]:\sH\to\left[\begin{array} {c} \sH'\\
H^2(\sD_{T'})
\end{array}\right].
\]
By virtue of $\tilde V\tilde BR_\circ=\tilde BQ_\circ$ it follows
that
\[
\left[\begin{array} {cc} T'&0\\
E_{\sD_{T'}}D_{T'}&S_{\sD_{T'}}
\end{array}\right]\left[\begin{array} {c} AR\\
\Gamma_1D_AR
\end{array}\right]=\left[\begin{array} {c} AQ\\
\Gamma_1D_AQ
\end{array}\right].
\]
Thus $B$ is a contraction, $A=\Pi_{\sH'}B$ and $VBR=BQ$, that is,
$B$ is a contractive interpolant for the data set
$\{A,T',V,R,Q\}$.

Next, let $B$ from $\sH$ to $\sV=\sH'\oplus H^2(\sD_{T'})$ be an
arbitrary contractive interpolant for the data set
$\{A,T',V,R,Q\}$. We have to show that $B$ is given by
(\ref{BBtil}), where $\tilde B$ is some contractive interpolant
for the data set (\ref{datatil1}). In fact, from (\ref{Btil1}) we
see that it suffices to find a contraction $\Gamma$ from $\sD_B$
into $H^2(\sD_\circ)$ such that the operator $\tilde B$, given by
\begin{equation}
\label{Btil2}  \tilde B=\left[\begin{array} {cc} B&0\\
0&I_{\sD_\circ}\\ \Gamma D_B&0\\
\end{array}\right] :
\left[\begin{array} {c} \sH\\\sD_{\circ}
\end{array}\right]\to
\left[\begin{array} {c} \sV\\\sD_{\circ}\\
H^2(\sD_\circ)
\end{array}\right],
\end{equation}
satisfies the intertwining relation $\tilde W\tilde
BR_\circ=\tilde BQ_\circ$, where $\tilde W$ is the operator which
one obtains by interchanging the second and the third column and
the second and third row in the operator matrix for $\tilde V$.
Put
\begin{equation}
\label{BVcirc} B_\circ=\left[\begin{array} {cc} B&0\\ 0&I_{\sD_\circ}\\
\end{array}\right]: \left[\begin{array} {c} \sH\\\sD_{\circ}
\end{array}\right]\to\left[\begin{array} {c} \sV \\\sD_{\circ}
\end{array}\right], \quad V_\circ=\left[\begin{array} {cc} V&0\\ 0&0
\end{array}\right] \mbox{ on } \left[\begin{array} {c} \sV\\\sD_{\circ}
\end{array}\right].
\end{equation}
Since $VBR=BQ$, we have $V_\circ B_\circ R_\circ=B_\circ Q_\circ$.
Now, notice that $B_\circ=B\oplus I_{\sD_\circ}$ is a contraction.
Furthermore, $V_\circ$  is a partial isometry, and the
Sz.-Nagy-Sch\"affer minimal isometric lifting of $V_\circ$ is
equal to $\tilde W$.  Thus
\begin{equation}
\label{datatil2} \big\{B_\circ, V_\circ,\tilde W,R_\circ, Q_\circ
\big\}
\end{equation}
is a lifting data set. Since $R_\circ^*R_\circ=Q_\circ^*Q_\circ$,
we know from Theorem \ref{mainth2} that the data set
(\ref{datatil2}) has a contractive interpolant $\tilde B$. By
identifying the spaces
\[\sH'\oplus\sD_\circ\oplus H^2(\sD_{T'}) \oplus
H^2(\sD_\circ)\quad\mbox{and}\quad \sH'\oplus
H^2(\sD_{T'})\oplus\sD_\circ \oplus H^2(\sD_\circ).
\]
one sees that  this operator $\tilde B$ is also a contractive
interpolant for the data set (\ref{datatil1}), and from
(\ref{Btil2}) it follows that with this choice of $\tilde B$ the
identity (\ref{BBtil}) holds. \epr

\medskip\noindent\textbf{Proof of Theorem \ref{mainth}.} We split
the proof into two parts.

\smallskip\noindent\textbf{Part 1.} Let $B$ be a contractive
interpolant for the data set $\{A,T',V,R,Q\}$. Then $B$ is of the
form (\ref{BBtil})  for some contractive interpolant $\tilde B$
for $\{A_\circ,T_\circ',\tilde V,R_\circ,Q_\circ\}$. Since
$R_\circ^*R_\circ=Q_\circ^*Q_\circ$, we can use Theorem
\ref{mainth2} to find  a formula for $\tilde B$. To write this
formula, we need the subspaces
\[
\sF_\circ=\overline{D_{A_\circ}Q_\circ\sH_0} \quad \mbox{and}\quad
\sF_\circ'=\overline{\left[\begin{array} {c} D_{T_\circ'} A_\circ
R_\circ\\ D_{A_\circ} R_\circ
\end{array}\right]\sH_0},
\]
and the unitary operator $\oo_\circ$ from $\sF_\circ$ onto
$\sF_\circ'$ given by
\[
\oo_\circ D_{A_\circ}Q_\circ=\left[\begin{array} {c}
D_{T_\circ'}A_\circ R_\circ\\ D_{A_\circ}R_\circ
\end{array}\right].
\]
In this setting,
\begin{equation}
\label{Dtil} D_{A_\circ}=\left[\begin{array} {cc} D_A&0\\0&0
\end{array}\right]
\quad \mbox{and}\quad D_{T_\circ'}=\left[\begin{array} {cc}
D_{T'}&0\\0&I_{\sD_\circ}
\end{array}\right].
\end{equation}
A straightforward computation shows that
\[
D_{A_\circ}Q_\circ=\left[\begin{array} {c} D_AQ \\ 0
\end{array}\right]:\sH_0\to \left[\begin{array} {c} \sH\\
\sD_\circ
\end{array}\right]
\]
and
\[
 \left[\begin{array} {c} D_{T_\circ'}A_\circ R_\circ\\
D_{A_\circ}R_\circ \end{array}\right]
=\left[\begin{array} {c} D_{T'}AR\\ D_\circ \\
D_AR\\0
\end{array}\right]:\sH_0\to \left[\begin{array} {c} \sH'\\ \sD_\circ\\
\sH\\ \sD_\circ \end{array}\right].
\]
By interchanging in the last column the first two coordinate
spaces and identifying the vector  $x\oplus0$ with the vector $x$,
we see that
\begin{equation}
\label{FFtil} \sF_\circ=\sF,\qquad \sF_\circ'=\sF' \quad
\mbox{and}\quad \oo_\circ=\oo,
\end{equation}
where the subspaces $\sF$ and $\sF'$ and the unitary operator
$\oo$ are defined in Section \ref{sec:intro}. Let us now apply
Theorem \ref{mainth2} to $\tilde B$. It follows that
\begin{equation}
\label{Btil3} \tilde Bx=\left[\begin{array} {c} A_\circ x\\
\Pi_{T_\circ'}F_\circ(\la)\big(I_{\sD_{A_\circ}}-\la\Pi_{A_\circ}F_\circ(\la)\big)^{-1}D_{A_\circ}x
\end{array}\right],\quad x\in \left[\begin{array} {c}
\sH\\ \sD_\circ \end{array}\right],
\end{equation}
where $F_\circ\in\eS(\sD_{A_\circ},\sD_{T_\circ'}\oplus
\sD_{A_\circ})$ satisfies $F_\circ(0)|\sF_\circ=\oo_\circ$. Here
\[
\Pi_{T_\circ'}=\left[\begin{array}{cc} I&0
\end{array}\right]: \left[
\begin{array}{c}
\sD_{T_\circ'}\\ \sD_{A_\circ}
\end{array}
\right]\to \sD_{T_\circ'}
 \quad\mbox{and}\quad
\Pi_{A_\circ}=\left[\begin{array}{cc} 0&I
\end{array}\right]: \left[
\begin{array}{c}
\sD_{T_\circ'}\\ \sD_{A_\circ}
\end{array}
\right]\to \sD_{A_\circ}.
\]
{}From (\ref{Dtil}) we see that we can identify in a canonical way
$\sD_{A_\circ}$ with $\sD_A$, and $\sD_{T_\circ'}$ with
$\sD_\circ\oplus\sD_{T'}$. This together with (\ref{FFtil}) shows
that we can view $F_\circ$ as a function $F$ from the Schur class
$\eS(\sD_A,\sD_\circ\oplus\sD_{T'}\oplus\sD_A)$ satisfying
$F(0)|\sF=\oo$ and
\begin{equation}
\label{Btil4} \tilde B \left[\begin{array} {c} h\\ d_0
\end{array}\right]=\left[\begin{array} {cc}
Ah&0\\ 0&d_0\\
\Pi_{T'}F(\la)(I-\la\Pi_AF(\la))^{-1}D_Ah&0\\
\Pi_{D_\circ}F(\la)(I-\la\Pi_AF(\la))^{-1}D_Ah&0
\end{array}\right],\quad \left[\begin{array} {c} h\\ d_0
\end{array}\right]\in\left[\begin{array} {c} \sH\\ \sD_\circ
\end{array}\right].
\end{equation}
Here $\Pi_{T'}$ and $\Pi_A$ are the projections given by
(\ref{projects}) and
\[
\Pi_{D_\circ}=\left[\begin{array} {ccc} I&0&0
\end{array}\right]:\left[\begin{array} {c} \sD_\circ\\ \sD_{T'}\\
\sD_A \end{array}\right]\to\sD_\circ.
\]
Since $B$ is obtained from $\tilde B$ via (\ref{BBtil}), we
conclude that $B$ has the desired form (\ref{sols}).

\smallskip\noindent\textbf{Part 2.} The reverse implication is
proved in a similar way. Indeed, assume that $B$ is given by
(\ref{sols}), where
$F\in\eS(\sD_A,\sD_\circ\oplus\sD_{T'}\oplus\sD_A)$ satisfies
$F(0)|\sF=\oo$. Using the identifications made in the first part
of the proof, we can  view  $F$ as a function
$F_\circ\in\eS(\sD_{A_\circ},\sD_{T_\circ'}\oplus\sD_{A_\circ})$
satisfying $F_\circ(0)|\sF_\circ=\oo_\circ$. But then we can use
Theorem \ref{mainth2} to show that $\tilde B$ defined by
(\ref{Btil3}) is a contractive interpolant for the data set
$\{A_\circ,T_\circ',\tilde V,R_\circ,Q_\circ\}$. Since $\tilde B$
is also given by (\ref{Btil4}), we conclude that $B$ and $\tilde
B$ are related as in (\ref{BBtil}). Thus Proposition
\ref{propbbtil} implies that $B$ is a contractive interpolant for
$\{A,T',V,R,Q\}$.\epr

\section{Parameterization and uniqueness of solutions}
\label{sec:unique} \setcounter{equation}{0}

\noindent In this section we prove the second main theorem
(Theorem~\ref{thmmain2}). As a consequence of this theorem we
obtain conditions on the lifting data set $\liftsetV$ guaranteeing
that the parameterization in Theorem~\ref{mainth} is proper, that
is, conditions on $\liftsetV$ implying that for every contractive
interpolant $B$ for $\liftsetV$ there exists a unique $F$ in
$\eS(\sD_A,\sD_\circ\oplus\sD_{T'}\oplus\sD_A)$ with
$F(0)|\sF=\oo$ such that $B=B_F$. Here $B_F$ is the contractive
interpolant for $\liftsetV$ produced by the Schur class function
$F$ from $\eS(\sD_A,\sD_\circ\oplus\sD_{T'}\oplus\sD_A)$ with
$F(0)|\sF=\oo$ as in Theorem~\ref{mainth}, that is,
\begin{equation}\label{eq:BF}
B_Fh=\mat{c}{ Ah\\ \Pi_{T'}F(\lambda)(I_{\sD_A}-\lambda \Pi_A
F(\lambda))^{-1}D_Ah}, \quad h\in \sH.
\end{equation}
We shall also present conditions on $\liftsetV$ implying the
existence  of a unique interpolant for $\liftsetV$.

To shorten the notation in this section we define
\begin{equation}
\sV=\sH'\oplus H^2(\sD_{T'})\quad\mbox{and}\quad
\tilde\sV=\sH'\oplus\sD_\circ\oplus H^2(\sD_{T'})\oplus
H^2(\sD_\circ).
\end{equation}
Also, for a given contractive interpolant $B$ for $\liftsetV$ we
define the spaces $\sF_B$ and $\sF_B'$ by
\begin{equation}\label{eq:FFB}
\sF_B=\overline{D_BQ\sH_0}\quad\mbox{and}\quad
\sF_B'=\overline{\mat{c}{D_\circ\\D_BR}\sH_0}.
\end{equation}
Notice that $\sG_B$ and $\sG_B'$ in~(\ref{eq:FB0}) are then given
by
\begin{equation}\label{eq:GB}
\sG_B=\sD_B\ominus\sF_B\quad\mbox{and}\quad
\sG_B'=(\sD_\circ\oplus\sD_B)\ominus\sF_B'.
\end{equation}
With the above notation and definitions we can reformulate
Theorem~\ref{thmmain2} as follows.

\begin{theorem}\label{th:main2}
Let $\liftsetV$ be a lifting data set with $V$ the
Sz.-Nagy-Sch\"affer minimal isometric lifting of $T'$, and let $B$
be a contractive interpolant for the data set $\liftsetV$. Then
there exists a one to one mapping from the set of all $F$ in
$\eS(\sD_A,\sD_\circ\oplus\sD_{T'}\oplus\sD_A)$ with
$F(0)|\sF=\oo$ such that $B=B_F$ onto the Schur class
$\eS(\sG_B,\sG_B')$, with $\sG_B$ and $\sG_B'$ as
in~$(\ref{eq:GB})$.
\end{theorem}

For the proof of Theorem~\ref{th:main2} it will be convenient to
first prove two lemma's. Let $\liftsetVT$ be as defined in
Section~\ref{sec:RCLT1}. Given  a contractive interpolant $B$ for
the data set $\liftsetV$, we let  $B_\circ$ and $V_\circ$ be the
operators defined by~(\ref{BVcirc}). Furthermore,  as in the
previous section, we define $\tilde W$ to be the operator which
one obtains by interchanging the second and the third column and
the second and third row in the operator matrix for $\tilde V$.
Recall that both $\liftsetVT$ and $\{B_\circ,V_\circ,\tilde
W,R_\circ,Q_\circ\}$ are lifting data sets. {}From the
construction of $\tilde W$ from $\tilde V$ we see that both
$\tilde W$ and $\tilde V$ are minimal isometric liftings of both
$T_\circ'$ and $V_\circ$.

\begin{lemma}\label{lem:VBtoTA}
Let $B$ be a contractive interpolant for $\liftsetV$ and let the
pair $\coupon{}$ be an isometric coupling of
$\{V_\circ,B_\circ\}$. Then
\begin{itemize}
\item[(i)] the pair $\coup{}$ is an isometric coupling of
$\lifsetT$;

\item[(ii)] the pair $\coup{}$ is minimal as an isometric coupling
of
            $\{V_\circ,B_\circ\}$ if and only if $\coup{}$ is minimal as an
            isometric coupling of $\lifsetT$;

\item[(iii)] the operator $U$ is an isometric lifting of both $T'$
and
             $T_\circ'$; moreover,
             $\Phi_{U,T'}$ is the canonical embedding of $\sV$ into $\sK$,
             and $\Phi_{U,T_\circ'}v=v$ for all $v\in\sV$;

\item[(iv)] the contractive interpolant
            $B=\Pi_\sV\Phi^*_{U,T_\circ'}\tau\Pi^*_\sH$.
\end{itemize}
Furthermore, two isometric couplings $\coup{1}$ and $\coup{2}$ of
$\{V_\circ,B_\circ\}$ are isomorphic as isometric couplings of
$\{V_\circ,B_\circ\}$ if and only if they are isomorphic as
isometric couplings of $\lifsetT$.
\end{lemma}

\bpr First remark that $\sK$ can be decomposed as
$\sV\oplus\sD_\circ\oplus\sM$ for some Hilbert space $\sM$.
Relative to this direct sum decomposition the operators $U$ and
$\tau$ admit operator matrix representations of the form
\begin{equation}\label{eq:Utaufin}
U=\mat{ccc}{V&0&0\\0&0&0\\0&*&*}\mbox{ on
}\mat{c}{\sV\\\sD_\circ\\\sM}\quad
\mbox{and}\quad\tau=\mat{cc}{B&0\\0&I_{\sD_\circ}\\ *&0}:
\mat{c}{\sH\\\sD_\circ}\to\mat{c}{\sV\\ \sD_\circ\\\sM},
\end{equation}
where $*$ represents operators which are not specified any
further.

\smallskip\noindent\textbf{(i)}
Because $V$ is an isometric lifting of $T'$ and $U$ is an
isometric lifting of $V$, as we can see from~(\ref{eq:Utaufin}),
we obtain that $U$ is an isometric lifting of $T'$.
{}From~(\ref{eq:Utaufin}) we can immediately see that $U$ also is
an isometric lifting of the zero operator on $\sD_\circ$. Hence
$U$ is an isometric lifting of $T_\circ'$. Since $\tau$ is as
in~(\ref{eq:Utaufin}) and $\Pi_{\sH'}B=A$, we see that
$\Pi_{\sH'\oplus\sD_\circ}\tau=A_\circ$. So $\coup{}$ is an
isometric coupling of $\lifsetT$.

\smallskip\noindent\textbf{(ii)}
Assume that $\coup{}$ is minimal as an isometric coupling of
$\{V_\circ,B_\circ\}$, and thus that the space
$(\sV\oplus\sD_\circ)\vee\tau(\sH\oplus\sD_\circ)$  is cyclic for
$U$.

Notice that in general we have for every operator $W$ on a Hilbert
space $\sL$ with $\sU$ and $\sY$ subspaces $\sL$ that
\[
\bvz W^n(\sU\vee\sY)=(\bvz W^n\sU)\vee(\bvz W^n\sY).
\]
Applying this with  $U$  in (\ref{eq:Utaufin}), the fact that  $V$
on $\sV$ is a minimal isometric lifting of $T'$ and the fact that
$(\sV\oplus\sD_\circ)\vee\tau(\sH\oplus\sD_\circ)$  is cyclic for
$U$ yield
\begin{eqnarray*}
\sK&=&\bvz U^n((\sV\oplus\sD_\circ)\vee\tau(\sH\oplus\sD_\circ))
\\
 &=&(\bvz U^n\sV)\vee(\bvz U^n(\sD_\circ\vee\tau(\sH\oplus\sD_\circ)))\\
&=&\sV\vee(\bvz U^n(\sD_\circ\vee\tau(\sH\oplus\sD_\circ)))\\
 &=&(\bvz\Pi^*_\sV V^n\sH')\vee(\bvz U^n(\sD_\circ\vee\tau(\sH\oplus\sD_\circ)))\\
&=&(\bvz U^n\sH')\vee(\bvz
U^n(\sD_\circ\vee\tau(\sH\oplus\sD_\circ)))\\ &=&\bvz
U^n((\sH'\oplus\sD_\circ)\vee\tau(\sH\oplus\sD_\circ)).
\end{eqnarray*}
Hence $(\sH'\oplus\sD_\circ)\vee\tau(\sH\oplus\sD_\circ)$ is
cyclic for $U$, and thus $\coup{}$ is minimal as an isometric
coupling of $\lifsetT$.

Conversely, assume that $\coup{}$ is minimal as an isometric
coupling of \linebreak $\lifsetT$. In other words,
$(\sH'\oplus\sD_\circ)\vee\tau(\sH\oplus\sD_\circ)$ is cyclic for
$U$. Note that $\sH'$ is a subspace of $\sV$ and hence
$(\sH'\oplus\sD_\circ)\vee\tau(\sH\oplus\sD_\circ)$ is a subspace
of $(\sV\oplus\sD_\circ)\vee\tau(\sH\oplus\sD_\circ)$. This
implies that $(\sV\oplus\sD_\circ)\vee\tau(\sH\oplus\sD_\circ)$ is
cyclic for $U$ as well. Thus $\coup{}$ is minimal as an isometric
coupling of $\{V_\circ,B_\circ\}$.

\smallskip\noindent\textbf{(iii)}
We already showed, in (i), that $U$ is an isometric lifting of
both $T'$ and $T_\circ'$. {}From~(\ref{eq:Utaufin}) we see that
$V$ is the minimal isometric lifting of $T'$ associated with $U$
and thus, using the remark in the final paragraph of
Section~\ref{sec:IL}, we obtain that
\[
\Phi_{U,T'}=\Pi^*_\sV\Phi_{V,T'}=\Pi^*_\sV I_\sV=\Pi^*_\sV.
\]

Because $\Phi_{U,T_\circ'}$ is the unique isometry associated with
$T_\circ'$ that intertwines $\tilde V$ with $U$, we see that the
isometry $\Phi_{U,T_\circ'}\Pi^*_\sV$  satisfies
$\Phi_{U,T_\circ'}\Pi^*_\sV h=h$ for all $h\in\sH'$ and
\[
U\Phi_{U,T_\circ'}\Pi^*_\sV=\Phi_{U,T_\circ'}\tilde V\Pi^*_\sV
=\Phi_{U,T_\circ'}\Pi^*_\sV V.
\]
Hence $\Phi_{U,T_\circ'}\Pi^*_\sV$ is the unique isometry
associated with $T'$ that intertwines $V$ and $U$. Thus for all
$v\in\sV$ we have
\[
\Phi_{U,T_\circ'}v=\Phi_{U,T_\circ'}\Pi^*_\sV v=\Phi_{U,T'}v=v.
\]

\smallskip\noindent\textbf{(iv)}
In the proof of (iii) we saw that
$\Phi_{U,T_\circ'}\Pi^*_\sV=\Phi_{U,T'}=\Pi^*_\sV$. Hence
from~(\ref{eq:Utaufin}) we obtain that
\[
B=\Pi_\sV\tau\Pi^*_\sH=\Pi_\sV\Phi^*_{U,T_\circ'}\tau\Pi^*_\sH.
\]

It remains to prove the final statement of the lemma. For this
purpose, let $\coupon{1}$ and $\coupon{2}$ be isometric couplings
of $\{V_\circ,B_\circ\}$.

First assume that $\coup{1}$ and $\coup{2}$ are isomorphic as
isometric couplings of $\{V_\circ,B_\circ\}$. We can immediately
see from the definition of an isomorphism and the fact that
$\sH'\oplus\sD_\circ$ is a subspace of $\sV\oplus\sD_\circ$, that
every isomorphism from $\coup{1}$ to $\coup{2}$ as isometric
couplings of $\{V_\circ,B_\circ\}$ also is an isomorphism from
$\coup{1}$ to $\coup{2}$ as isometric couplings of $\lifsetT$.
Hence $\coup{1}$ and $\coup{2}$ are isomorphic as isometric
couplings of $\lifsetT$.

Conversely, assume that $\coup{1}$ and $\coup{2}$ are isomorphic
as isometric couplings of $\lifsetT$ and that $\Psi$  is an
isomorphism from $\coup{1}$ to $\coup{2}$. Then $\Psi\Pi^*_\sV$ is
an isometry from $\sV$ to $\sK_2$ with $\Psi\Pi^*_\sV h=h$ for
each $h\in\sH'$. Since $V$ is the minimal isometric lifting of
$T'$ associated with $U_1$, we obtain
\[
U_2\Psi\Pi^*_\sV=\Psi U_1\Pi^*_\sV=\Psi\Pi^*_\sV V.
\]
Thus with (iii) we see that
\[
\Psi\Pi^*_\sV=\Phi_{U_2,T'}=\Pi^*_\sV.
\]
Since $\Psi$ is an isomorphism between isometric couplings of
$\lifsetT$, the operator $\Psi$ is the identity on
$\sH'\oplus\sD_\circ$. In particular, $\Psi d=d$ for each
$d\in\sD_\circ$. Hence $\Psi$ also is an isomorphism from
$\coup{1}$ to $\coup{2}$ as isometric couplings of
$\{V_\circ,B_\circ\}$. \epr

\begin{lemma}\label{lem:TAtoVB}
Let $B$ be a contractive interpolant for the data set $\liftsetV$,
and let $\coupon{}$ be an isometric coupling of $\lifsetT$ such
that
\[
B=\Pi_\sV\Phi^*_{U,T_\circ'}\tau\Pi^*_\sH.
\]
Then there exists an
isometric coupling $\{\check U,\check\tau\}$ of $\lifsetT$,
isomorphic to $\coup{}$, such that $\{\check U,\check\tau\}$ also
is an isometric coupling of $\{V_\circ,B_\circ\}$.
\end{lemma}

\bpr {}From the remark in the last paragraph of
Section~\ref{sec:IL} we can conclude that
$\Phi_{U,T_\circ'}=\Pi^*_{\sK'}\Phi_{U',T_\circ'}$, where $U'$ on
$\sK'$ is the minimal isometric lifting of $T_\circ'$ associated
with $U$. Since $U'$ is minimal and $\tilde V$ is the
Sz.-Nagy-Sch\"affer minimal isometric lifting of $T_\circ'$, we
see that the unique isometry $\Phi_{U',T_\circ'}$ associated with
$T_\circ'$ that intertwines $\tilde V$ with $U'$ is unitary.
Define $\check\Psi$ from $\tilde\sV\oplus\sW$ to $\sK$, where
$\sW=\sK\ominus\sK'$, by
\[
\check\Psi=\mat{cc}{\Phi_{U',T_\circ'}&0\\0&I_\sW}:\mat{c}{\tilde\sV\\\sW}\to
\mat{c}{\sK'\\\sW}.
\]
Then $\check\Psi$ is a unitary operator with $\check\Psi x=x$ for
all $x\in\sH'\oplus\sD_\circ$. So the operators $\check
U=\check\Psi^*U\check\Psi$ and $\check\tau=\check\Psi^*\tau$ form
an isometric coupling $\{\check U,\check\tau\}$ of $\lifsetT$ that
is isomorphic to $\coup{}$.

Since $\Phi_{U',T_\circ'}$ intertwines $\tilde V$  and $U'$,  we
obtain that $\tilde V$ is the minimal isometric lifting of
$T_\circ'$ associated with $\check U$. Recall that $\tilde V$ also
is an isometric lifting of $V_\circ$. Hence $\check U$ is an
isometric lifting of $V_\circ$.

Because $\tilde V$ is the minimal isometric lifting of $T_\circ'$
associated with $\check U$, we have  $\Phi_{\check
U,T_\circ'}=\Pi^*_{\tilde\sV}$. Since   $\{\check U,\check\tau\}$
and $\coup{}$ are isomorphic, Lemma \ref{lemisocoup12} implies
that
\[
B=\Pi_\sV\Phi^*_{U,T_\circ'}\tau\Pi^*_\sH
 =\Pi_\sV\Phi^*_{\check U,T_\circ'}\check\tau\Pi^*_\sH
 =\Pi_\sV\Pi_{\tilde\sV}\check\tau\Pi^*_\sH
 =\Pi_\sV\check\tau\Pi^*_\sH.
\]
Note that $\Pi_{\sH'\oplus\sD_\circ}\check\tau=A_\circ$, and thus,
since both $\check\tau$ and $A_\circ|\sD_\circ$ are isometries, we
get that $\Pi_{\sV\oplus\sD_\circ}\check\tau=B_\circ$. Hence
$\{\check U,\check\tau\}$ is an isometric coupling of
$\{V_\circ,B_\circ\}$. \epr

\medskip\noindent\textbf{Proof of Theorem \ref{th:main2}.} Let
$\eS_B$ be the set defined by
\begin{equation}\label{eq:SB}
\eS_B=\{ F\in\eS(\sD_A,\sD_\circ\oplus\sD_{T'}\oplus\sD_A) \mid
F(0)|\sF=\oo \mbox{ and }B=B_F \}.
\end{equation}
We have to show that there exists a one to one mapping from
$\eS_B$ onto $\eS(\sG_B,\sG_B')$.

By applying Theorem~\ref{thm:ic-schur} to the pair $\lifsetT$ and
Proposition~\ref{ic-schurRQ} to the lifting data set $\liftsetVT$,
and using the identities in~(\ref{FFtil}), we obtain that the
mapping
\begin{equation} \label{mapUtauF}
\coup{}\mapsto F_{\coup{}}
\end{equation}
given by Theorem~\ref{thm:ic-schur} is a one to one mapping from
the set of (equivalence classes of) minimal isometric couplings
$\coup{}$ of $\lifsetT$ satisfying $U\tau R_\circ=\tau Q_\circ$
onto the set of all functions
$F\in\eS(\sD_A,\sD_\circ\oplus\sD_{T'}\oplus\sD_A)$ satisfying
$F(0)|\sF=\oo$. Moreover, from Proposition~\ref{prop:couplsol} and
Proposition~\ref{propbbtil}, applied to $\liftsetVT$, we obtain
that the mapping (\ref{mapUtauF}) maps the set of (equivalence
classes of) minimal isometric couplings $\coup{}$ of $\lifsetT$
satisfying $U\tau R_\circ=\tau Q_\circ$ and
$B=\Pi_\sV\Phi^*_{U,T_\circ'}\tau\Pi^*_\sH$ onto the set $\eS_B$
defined by (\ref{eq:SB}).

Then, using Lemma~\ref{lem:VBtoTA} and Lemma~\ref{lem:TAtoVB}, we
obtain that there exists a one to one mapping from $\eS_B$ onto
the set of (equivalence classes of) minimal isometric couplings
$\coup{}$ of $\{V_\circ,B_\circ\}$ satisfying $U\tau R_\circ=\tau
Q_\circ$.

Note that, because $B$ is a contractive interpolant for
$\liftsetV$ and thus $VBR=BQ$, we have that $\{B,V,V,R,Q\}$ is a
lifting data set, and that the lifting data set
$\{B_\circ,V_\circ,\tilde W,R_\circ,Q_\circ\}$ is constructed from
$\{B,V,V,R,Q\}$ in the same way as we constructed $\liftsetVT$
from $\liftsetV$ in Section~\ref{sec:RCLT1}. Moreover, since $V$
is an isometry and thus $\sD_V=\{0\}$, we get that $\sF_B$ and
$\sF_B'$ in~(\ref{eq:FFB}) correspond to $\{B,V,V,R,Q\}$ as $\sF$
and $\sF'$ correspond to $\liftsetV$. Hence there exists a unique
unitary operator $\oo_B$ from $\sF_B$ to $\sF_B'$ defined by
\[
\oo_BD_BQ=\mat{c}{D_\circ\\D_BR}.
\]

By again applying Theorem~\ref{thm:ic-schur},
Proposition~\ref{ic-schurRQ} and the identities in~(\ref{FFtil}),
but now to the pair $\{V_\circ,B_\circ\}$ and the lifting data set
$\{B_\circ,V_\circ,\tilde W,R_\circ,Q_\circ\}$, we obtain that
there exists a one to one mapping from the set of (equivalence
classes of) minimal isometric couplings $\coup{}$ of
$\{B_\circ,V_\circ,\tilde W,R_\circ,Q_\circ\}$ satisfying $U\tau
R_\circ=\tau Q_\circ$ onto the set of all functions
$H\in\eS(\sD_B,\sD_\circ\oplus\sD_B)$ satisfying
$H(0)|\sF_B=\oo_B$. Thus there exists a one to one mapping from
$\eS_B$ onto the set of all functions
$H\in\eS(\sD_B,\sD_\circ\oplus\sD_B)$ satisfying
$H(0)|\sF_B=\oo_B$.

For each $H\in\eS(\sD_B,\sD_\circ\oplus\sD_B)$ we have that
$H(0)|\sF_B=\oo_B$ if and only if there exists a (unique) $G$ in
$\eS(\sG_B,\sG_B')$, with $\sG_B$ and $\sG_B'$ as
in~(\ref{eq:GB}), such that
\begin{equation}\label{eq:HtoG}
H(\la)=\mat{cc}{\oo_B&0\\0&G(\la)}:
\mat{c}{\sF_B\\\sG_B}\to\mat{c}{\sF_B'\\\sG_B'},\quad\la\in\mathbb{D}.
\end{equation}
Hence there exists a one to one mapping from the set $\eS_B$ onto
$\eS(\sG_B,\sG_B')$. \epr

\medskip
In fact, in the proof of Theorem~\ref{th:main2} we do not
only show that there exists a one to one mapping from the set of
$F$ in $\eS(\sD_A,\sD_\circ\oplus\sD_{T'}\oplus\sD_A)$ with
$F(0)|\sF=\oo$ such that $B=B_F$ onto $\eS(\sG_B,\sG_B')$, but we
actually indicate how such a mapping can be constructed. To be
more specific, the construction in the reverse way goes as
follows.

Assume that $G$ is a Schur class function from
$\eS(\sG_B,\sG_B')$, with B some contractive interpolant for
$\liftsetV$. Define $H\in\eS(\sD_B,\sD_\circ\oplus\sD_B)$
by~(\ref{eq:HtoG}). Then $H$ satisfies $H(0)|\sF_B=\oo_B$, and
thus from Section~\ref{sec:IR} we obtain that there exists an
isometry $M$ from $\sD_\circ\oplus\sD_B\oplus\sY$ to
$\sD_B\oplus\sY$, for some Hilbert space $\sY$, such that
\[
H(\la)=\Pi_{\sD_\circ\oplus\sD_B}M^*(I_{\sD_B\oplus\sY}-\la
J_\sY^*M^*)^{-1} \Pi_{\sD_B}^*, \quad\la\in\mathbb{D},
\]
where $M$  satisfies the controllability type condition
\begin{equation}\label{eq:Mmin}
\sY=\Pi_\sY\bvz(J_\sY
M)^n\mat{c}{\sD_\circ\oplus\sD_B\\\{0\}}\quad\mbox{and} \quad
M|\sF_B'=\oo_B^*.
\end{equation}
Here $J_\sY$ is the partial isometry given by
\[
J_\sY=\mat{cc}{0&0\\0&I_\sY}:
\mat{c}{\sD_B\\\sY}\to\mat{c}{\sD_\circ\oplus\sD_B\\\sY}.
\]
Notice that because $V_\circ$ and $B_\circ$ are as
in~(\ref{BVcirc}) we obtain that $D_{V_\circ}=\Pi_{\sD_\circ}$,
$\sD_{V_\circ}=\sD_\circ$, $D_{B_\circ}=D_B$ and
$\sD_{B_\circ}=\sD_B\Pi_\sH$. Thus we can define
\[
\check
U=\mat{cc}{V_\circ&0\\M|\sD_\circ\Pi_{\sD_\circ}&M|(\sD_B\oplus\sY)}
\mbox{ on }\mat{c}{\sV\oplus\sD_\circ\\\sD_B\oplus\sY}
\]
and
\[
\check\tau=\mat{c}{B_\circ\\\Pi^*_{\sD_B}D_B\Pi_\sH}:\sH\oplus\sD_\circ\to
\mat{c}{\sV\oplus\sD_\circ\\\sD_B\oplus\sY}.
\]
Then $\{\check U,\check\tau\}$ is a  special isometric coupling of
$\{V_\circ,B_\circ\}$.  Because $M$ satisfies~(\ref{eq:Mmin}), the
coupling $\{\check U,\check\tau\}$ is minimal and $\check
U\check\tau R_\circ=\check\tau Q_\circ$. Hence by
Lemma~\ref{lem:VBtoTA} we obtain that $\{\check U,\check\tau\}$
also is a minimal isometric coupling of $\lifsetT$ with
$B=\Pi_\sV\Phi_{\check U,T_\circ'}\check\tau\Pi_\sH^*$. According
to Proposition~\ref{prop:special}, the coupling $\{\check
U,\check\tau\}$ is isomorphic to a special isometric coupling
$\coupon{}$ of $\lifsetT$. This isometric coupling $\coup{}$ is
minimal, satisfies $U\tau R_\circ=\tau Q_\circ$ and, by
Lemma~\ref{lemisocoup12}, we have that $B =
\Pi_\sV\Phi^*_{U,T_\circ'}\Pi_\sH^*$. The isometry $U$ defines an
isometry $Y$ from $\sD_{T'}\oplus\sD_\circ\oplus\sD_A\oplus\sX$ to
$\sD_A\oplus\sX$, with
$\sX=\sK\ominus(\sH'\oplus\sD_\circ\oplus\sD_A)$,
by~(\ref{eq:U1}), with $T_\circ'$ on $\sH'\oplus\sD_\circ$ instead
of $T'$ on $\sH'$. Then the functions $F$ in
$\eS(\sD_A,\sD_\circ\oplus\sD_{T'}\oplus\sD_A)$ with
$F(0)|\sF=\oo$ satisfying $B=B_F$ corresponding to the function
$G\in\eS(\sG_B,\sG_B')$ is given by
\[
F(\la)=\Pi_{\sD_{T'}\oplus\sD_\circ\oplus\sD_A}Y^*(I-\la
J_\sX^*Y^*)^{-1} \Pi_{\sD_A},\quad\la\in\mathbb{D},
\]
with $J_\sX$ being the partial isometry given by
\[
J_\sX=\mat{cc}{0&0\\0&I_\sX}:
\mat{c}{\sD_A\\\sX}\to\mat{c}{\sD_{T'}\oplus\sD_\circ\oplus\sD_A\\\sX}.
\]

{}From Theorem~\ref{th:main2} we immediately obtain the next
corollary.

\begin{corollary}\label{cor:uniq}
Let $B$ be a contractive interpolant for $\liftsetV$. Then there
is a unique $F$ in $\eS(\sD_A,\sD_\circ\oplus\sD_{T'}\oplus\sD_A)$
with $F(0)|\sF=\oo$ such that $B=B_F$ if and only if $\sF_B=\sD_B$
or $\sF_B'=\sD_\circ\oplus\sD_B$.
\end{corollary}

The next lemma gives some sufficient conditions on $\liftsetV$
under which the parameterization in Theorem~\ref{mainth} is
proper. To this end, define the subspace $\sF_A'$ of
$\sD_\circ\oplus\sD_A$ by
\begin{equation}
\sF_A'=\overline{\mat{c}{D_\circ\\D_AR}\sH_0}.
\end{equation}

\begin{lemma}\label{lem:prop}
Let $\liftsetV$ be a lifting data set. Then
\begin{itemize}
\item[(i)] $\sF=\sD_A$ implies that there exists a unique
contractive
           interpolant $B$ and that $\sF_B=\sD_B$;

\item[(ii)] $\sF_A'=\sD_\circ\oplus\sD_A$ implies that
            $\sF_B'=\sD_\circ\oplus\sD_B$ for every contractive interpolant $B$.
\end{itemize}
If either $\sF=\sD_A$ or $\sF_A'=\sD_\circ\oplus\sD_A$ holds, then
the mapping $F\rightarrow B_F$ given by~$(\ref{eq:BF})$ is one to
one from the set of all $F\in  \eS(\sD_A,
\sD_\circ\oplus\sD_{T'}\oplus \sD_A)$ satisfying $F(0)|\sF=\omega$
onto the set of all contractive interpolants.
\end{lemma}

\bpr Let $B$ be a contractive interpolant for $\liftsetV$. Then
$\Pi_{\sH'}B=A$, hence there exists a contraction $\Gamma$ from
$\sD_A$ to
 $H^2(\sD_{T'})$ such that
\[
B=\mat{c}{A\\\Gamma D_A}:\sH\to\mat{c}{\sH'\\ H^2(\sD_{T'})}.
\]
{}From this we obtain that for all $h\in\sH$
\begin{eqnarray*}
\|D_Bh\|^2 &=&\|h\|^2-\|Bh\|^2=\|h\|^2-\|Ah\|^2-\|\Gamma
D_Ah\|^2\\
&=&\|D_Ah\|^2-\|\Gamma D_A\|^2=\|D_\Gamma D_Ah\|^2.
\end{eqnarray*}
Note that $\sD_\Gamma\subset\sD_A$ and thus $\overline{D_\Gamma
D_A\sH}=\sD_\Gamma$. Hence there exists a unitary operator
$\gamma$ from $\sD_\Gamma$ onto $\sD_B$ such that $D_B=\gamma
D_\Gamma D_A$.

\smallskip\noindent\textbf{(i)}
Assume that $\sF=\sD_A$. Then there is only one $F$ in
$\eS(\sD_A,\sD_\circ\oplus\sD_{T'}\oplus\sD_A)$ with
$F(0)|\sF=\oo$. Hence, by Theorem~\ref{mainth}, there can be only
one contractive interpolant $B$ for $\liftsetV$, and for this
contractive interpolant $B$ there can be only one $F$ in
$\eS(\sD_A,\sD_\circ\oplus\sD_{T'}\oplus\sD_A)$ with
$F(0)|\sF=\oo$ such that $B=B_F$. Moreover, we have
\[
\sF_B =\overline{D_BQ\sH_0} =\overline{\gamma D_\Gamma D_AQ\sH_0}
=\overline{\gamma D_\Gamma\overline{D_AQ\sH_0}} =\overline{\gamma
D_\Gamma\sF} =\overline{\gamma D_\Gamma\sD_A} =\sD_B.
\]

\smallskip\noindent\textbf{(ii)}
Assume that $\sF_A'=\sD_\circ\oplus\sD_A$. Then we have that
\begin{eqnarray*}
\sF_B' &=&\overline{\mat{c}{D_\circ\\D_BR}\sH_0}
=\overline{\mat{cc}{I_{\sD_\circ}&0\\0&\gamma
D_\Gamma}\mat{c}{D_\circ\\ D_AR}\sH_0}\\
&=&\overline{\mat{cc}{I_{\sD_\circ}&0\\0&\gamma D_\Gamma}
\overline{\mat{c}{D_\circ\\ D_AR}\sH_0}}=
\overline{\mat{cc}{I_{\sD_\circ}&0\\0&\gamma D_\Gamma}\sF_A'}\\
&&\\
&=& \overline{\mat{cc}{I_{\sD_\circ}&0\\0&\gamma
D_\Gamma}\mat{c}{\sD_\circ\\\sD_A}}
=\mat{c}{\sD_\circ\\\overline{\gamma D_\Gamma \sD_A}}
=\mat{c}{\sD_\circ\\\sD_B}.
\end{eqnarray*}

\smallskip
The final statement of the lemma follows immediately from (i),
(ii) and Corollary~\ref{cor:uniq}.\epr

\medskip
For the classical commutant lifting theorem, that is, when
$\sH=\sH_0$, $R=I_\sH$ and $Q$ is an isometry on $\sH$, we have
already seen in Section~\ref{sec:intro} that the parameterization
in Theorem~\ref{mainth} is proper. This result also follows from
Lemma~\ref{lem:prop} (ii). Indeed, if $\sH=\sH_0$, $R=I_\sH$ and
$Q$ is an isometry on $\sH$, then
\[
\sF_A'=\overline{\{0\}\oplus D_A I_\sH \sH}
=\{0\}\oplus\sD_A=\sD_\circ\oplus\sD_A.
\]

Finally we derive some sufficient conditions on $\liftsetV$
guaranteeing that there is only one contractive interpolant.
{}From Lemma~\ref{lem:prop} we already know that the condition
$\sF=\sD_A$ is such a condition. In the same way we can see that
the condition $\sF'=\sD_{T'}\oplus\sD_\circ\oplus\sD_A$ is
sufficient.

For the classical commutant lifting theorem the combination of
these two conditions is also a necessary condition. That is, if
$\sH=\sH_0$, $R=I_\sH$ and $Q$ is an isometry on $\sH$, then there
is only one contractive interpolant if and only if $\sF=\sD_A$ or
$\sF'=\sD_{T'}\oplus\sD_\circ\oplus\sD_A$. We can see this as
follows. If the parametrization in Theorem~\ref{mainth} for the
lifting data set $\liftsetV$ is proper, then there is only one
contractive interpolant for $\liftsetV$ if and only if there is
only one $F$ in $\eS(\sD_A,\sD_\circ\oplus\sD_{T'}\oplus\sD_A)$
with $F(0)|\sF=\oo$. The latter is equivalent to the condition
`$\sF=\sD_A$ or $\sF'=\sD_{T'}\oplus\sD_\circ\oplus\sD_A$'.

Notice that when $T'$ is an isometry, then the Sz.-Nagy-Sch\"affer
minimal isometric lifting of $T'$ is $T'$ itself. So in that case
there also is only one contractive interpolant $B$ for
$\liftsetV$, namely $B=A$. In the next lemma we summarize the
above, and improve the condition
$\sF'=\sD_{T'}\oplus\sD_\circ\oplus\sD_A$ a bit further.

\begin{proposition}
Assume that for $\liftsetV$ either $T'$ is an isometry,
$\sF=\sD_A$ or $\sD_{T'}\oplus\sD_A\subset\sF'$. Then there exists
a unique contractive interpolant for $\liftsetV$.
\end{proposition}

\bpr We have already seen above that the requirement $T'$ is an
isometry and the equality $\sF=\sD_A$ are both sufficient
conditions. So assume that we have
$\sD_{T'}\oplus\sD_A\subset\sF'$. Define for all
$F\in\eS(\sD_A,\sD_\circ\oplus\sD_{T'}\oplus\sD_A)$ the Schur
class functions $F_\circ=\Pi_{\sD_\circ}F$,
$F_{T'}=\Pi_{\sD_{T'}}F$ and $F_A=\Pi_{\sD_A}F$. Hence for all
$\la\in\mathbb{D}$
\[
F(\la)=\mat{c}{F_\circ(\la)\\F_{T'}(\la)\\F_A(\la)}:\sD_A\to
\mat{c}{\sD_\circ\\\sD_{T'}\\\sD_A}.
\]
Then we have
\begin{equation}\label{eq:F0T'A}
\Pi_{T'}F(\la)(I_{\sD_A}-\la\Pi_AF(\la))^{-1}D_A
=F_{T'}(\la)(I_{\sD_A}-\la
F_A(\la))^{-1}D_A,\quad\la\in\mathbb{D}.
\end{equation}
All $F\in\eS(\sD_A,\sD_\circ\oplus\sD_{T'}\oplus\sD_A)$ with
$F(0)|\sF=\oo$ admit a matrix representation of the form
\[
F(\la)=\mat{cc}{\oo&0\\0&G(\la)}:\mat{c}{\sF\\\sG}\to\mat{c}{\sF'\\\sG'},\quad
\la\in\mathbb{D},
\]
for some $G\in\eS(\sG,\sG')$ where $\sG=\sD_A\ominus\sF$ and
$\sG'=(\sD_\circ\oplus\sD_{T'}\oplus\sD_A)\ominus\sF'$. Hence,
because $\sD_{T'}\oplus\sD_A\subset\sF'$ all
$F\in\eS(\sD_A,\sD_\circ\oplus\sD_{T'}\oplus\sD_A)$ with
$F(0)|\sF=\oo$ have identical $F_{T'}$ and $F_A$ and thus
from~(\ref{eq:F0T'A}) we see that $B_F$ is the same operator for
all $F\in\eS(\sD_A,\sD_\circ\oplus\sD_{T'}\oplus\sD_A)$ with
$F(0)|\sF=\oo$. Hence by Theorem~\ref{mainth} there is only one
contractive interpolant. \epr

\smallskip\noindent\textbf{Acknowledgement.} We thank Ciprian Foias for
useful discussions on an earlier version of this paper.

\smallskip\noindent\textbf{Added in proof.} At the IWOTA 2004
conference in Newcastle, when a preliminary version of this
paper had been completed, the authors
learned that W.S. Li and D. Timotin had a preprint ready in which
the coupling method was also used to study the relaxed commutant
lifting problem and the set of its solutions. Although the same method
was used in the same area the two papers turned out to be quite
complementary in style and results. We are happy that the editors
of \emph{Integral Equations and Operator Theory} agreed to publish the
final versions of both papers, one directly after the other in this
volume.


\end{document}